\documentclass[anglais]{cedram-aif}

\usepackage{version}
\excludeversion{annalif}
\includeversion{arxiv}

\usepackage[english]{babel}
\usepackage[latin1]{inputenc}
\usepackage[T1]{fontenc}

\usepackage{amssymb,amsmath,amsthm}
\usepackage{epic,eepic,longtable}
\usepackage{dsfont,mathrsfs}
\usepackage[neveradjust]{paralist}
\usepackage{url}
\input{xypic}
\xyoption{all}

\newcommand{\bu}{\bullet}
\newcommand{\ssm}{\smallsetminus}
\newcommand{\beq}{\begin{equation}}
\newcommand{\eeq}{\end{equation}}
\newcommand{\eps}{\varepsilon}
\newcommand{\ideal}[1]{{\langle#1\rangle}}
\newcommand{\into}{\hookrightarrow}
\newcommand{\onto}{\twoheadrightarrow}
\newcommand{\w}{\wedge}
\newcommand{\p}{\partial}
\newcommand{\xymat}{\SelectTips{cm}{}\xymatrix}

\newcommand{\del}[2]{\frac{\p}{\p x_{#1,#2}}}

\newcommand{\vp}{\varphi}
\newcommand{\Ga}{\Gamma}
\newcommand{\Om}{\Omega}
\newcommand{\om}{\omega}

\renewcommand{\mathbb}{\mathds}
\newcommand{\QQ}{\mathbb{Q}}
\newcommand{\RR}{\mathbb{R}}
\newcommand{\CC}{\mathbb{C}}

\renewcommand{\mathcal}{\mathscr}

\newcommand{\D}{\mathcal{D}}
\renewcommand{\O}{\mathcal{O}}

\newcommand{\gl}{\mathfrak{gl}}

\renewcommand{\sl}{\mathfrak{sl}}

\newcommand{\bb}{\mathfrak{b}}

\renewcommand{\gg}{\mathfrak{g}}

\renewcommand{\L}{\mathrm{L}}
\newcommand{\mm}{\mathfrak{m}}

\newcommand{\ttt}{\mathfrak{t}}

\renewcommand{\d}{\mathbf{d}}

\DeclareMathOperator{\ad}{ad}

\DeclareMathOperator{\Der}{Der}
\DeclareMathOperator{\diag}{diag}

\DeclareMathOperator{\Gl}{Gl}
\DeclareMathOperator{\gr}{gr}
\DeclareMathOperator{\Hom}{Hom}

\DeclareMathOperator{\Rep}{Rep}
\DeclareMathOperator{\sgn}{sign}
\DeclareMathOperator{\tr}{trace}
\DeclareMathOperator{\Sing}{Sing}
\DeclareMathOperator{\Sl}{Sl}

\DeclareMathOperator{\Sym}{Sym}

\DeclareMathOperator{\wt}{wt}
\DeclareMathOperator{\Z}{Z}

\title[Linear free divisors]{Linear free divisors and the\\global logarithmic comparison theorem}
\alttitle{Diviseurs lin\'eairement libres\\ et le th\'eor\`eme de\\ comparaison logarithmique global}

\author{\firstname{Michel} \lastname{Granger}}
\address{
Departement de Math\'ematiques\\
Universit\'e d'Angers\\
2 Bd Lavoisier\\
49045 Angers\\
France
}
\email{granger@univ-angers.fr}
\thanks{}

\author{\firstname{David} \lastname{Mond}}
\address{
Mathematics Institute\\
University of Warwick\\
Coventry CV47AL\\
England
}
\email{D.M.Q.Mond@warwick.ac.uk}
\thanks{DM is grateful to Ignacio de Gregorio for helpful conversations on the topics treated here.}

\author{\firstname{Alicia} \lastname{Nieto-Reyes}}
\address{
Departamento de Matematicas,\\
Estadistica y Computacion\\
Universidad de Cantabria\\
Spain
}
\email{alicia.nieto@unican.es}
\thanks{}

\author{\firstname{Mathias} \lastname{Schulze}}
\address{
Department of Mathematics\\
Oklahoma State University\\
Stillwater, OK 74078\\
United States}
\email{mschulze@math.okstate.edu}
\thanks{MS gratefully acknowledges financial support from EGIDE and the Humboldt Foundation.}

\thanks{We are grateful to the referee for a very careful reading and many valuable suggestions.}

\altkeywords{diviseur lin\'eairement libre, espace vectorielle pr\'ehomog\`ene, cohomologie de de~Rham, th\'eor\`eme de comparaison logarithmique, cohomologie des alg\`ebres de Lie, repr\'esentation des quivers}


\date{January 16, 2008}

\begin{document}

\begin{abstract}
A complex hypersurface $D$ in $\CC^n$ is a \emph{linear free divisor} (LFD) if its module of logarithmic vector fields has a global basis of linear vector fields. We classify all LFDs for $n$ at most $4$. 

By analogy with Grothendieck's comparison theorem, we say that the \emph{global logarithmic comparison theorem} (GLCT) holds for $D$ if the complex of global logarithmic differential forms computes the complex cohomology of $\CC^n\setminus D$.
We develop a general criterion for the GLCT for LFDs and prove that it is fulfilled whenever the Lie algebra of linear logarithmic vector fields is reductive.
For $n$ at most $4$, we show that the GLCT holds for all LFDs. 

We show that LFDs arising naturally as discriminants in quiver representation spaces (of real Schur roots) fulfill the GLCT. As a by-product we obtain a topological proof of a theorem of V.~Kac on the number of irreducible components of such discriminants.
\end{abstract}

\begin{altabstract}
Une hypersurface complexe de $\CC^n$ est appel\'ee un {\it diviseur lin\'eairement libre} (ou DLL) si son module de champs de vecteur logarithmiques a une base globale form\'ee de champs de vecteurs lin\'eaires. 
Nous classifions tous les DLL pour $n$ au plus \'egal a $4$.

Par analogie avec le th\'eor\`eme de comparaison de Grothendieck, on dit que le {\it th\'eor\`eme de comparaison logarithmique global} (ou TCLG) est vrai pour $D$ si le complexe des formes diff\'erentielles logarithmiques globales permet de calculer la 
cohomologie de $\CC^n\setminus D$ \`a coefficients complexes. 
Nous mettons en \'evidence un crit\`ere g\'en\'eral pour qu'un DLL ait la propri\'et\'e  TCLG, et nous d\'emontrons que ce crit\`ere s'applique lorsque l'alg\`ebre de Lie des champs de vecteurs logarithmiques lin\'eaires est r\'eductive. 
Pour $n$ inf\'erieur ou \'egal \`a $4$, nous montrons que le TCLG est vrai pour tous les DLL.

Nous montrons que les DLL qui apparaissent naturellement comme discriminants dans les espaces de repr\'esentations de carquois pour des racines de Schur r\'eelles satisfont au TCLG. 
Comme corollaire nous obtenons une d\'emonstration topologique d'un r\'esultat de V. Kac sur le nombre de composantes irr\'eductibles de tels discriminants.
\end{altabstract}

\maketitle

\section{Introduction}

We denote by $\O=\O_{\CC^n}$ the sheaf of holomorphic functions on $\CC^n$, by $\mm_p\subseteq\O_p$ the maximal ideal at $p\in\CC^n$, by $\Der=\Der_{\CC}(\O)$ the sheaf of $\CC$-linear derivations of $\O$ (or \emph{holomorphic vector fields}) on $\CC^n$, and by $\Om^\bullet=\Om^\bullet_{\CC^n}$ the complex of sheaves of holomorphic differential forms.
We shall frequently use a local or global coordinate system $x=x_1,\dots,x_n$ on $\CC^n$ and then denote by $\p=\p_1,\dots,\p_n$ the corresponding operators of partial derivatives $\p_i=\frac{\p}{\p x_i}$, $i=1,\dots,n$.
Note that $\Der=\bigoplus\O\cdot\p_i$ is a free $\O$-module of rank $n$.

Let $D\subseteq\CC^n$ be a reduced divisor.
K.~Saito \cite{Sai80} associated to $D$ the (coherent) sheaf of \emph{logarithmic vector fields} $\Der(-\log D)\subseteq\Der$ and the complex of (coherent) sheaves $\Om^\bu(\log D)\subseteq\Om^\bu(*D)$ of \emph{logarithmic differential forms} along $D$.
For a (local or global) defining equation $\Delta\in\O$ of the germ $D$, $\delta\in\Der$ is in $\Der(-\log D)$ if $\delta(\Delta)\in\O\cdot\Delta$, and $\om\in\Om^\bullet[\Delta^{-1}]$ is in $\Om^\bu(\log D)$ if $\Delta\cdot\om,\Delta\cdot d\om\in\Om^\bullet$.
Note that $\Der(-\log D)$ contains the \emph{annihilator} $\Der(-\log\Delta)$ of $\Delta$ defined by the condition $\delta(\Delta)=0$.
Saito showed that $\Der(-\log D)$ and $\Om^1(\log D)$ are reflexive and mutually dual and introduced the following important class of divisors.

\begin{defi}
A divisor $D$ is called \emph{free} if $\Der(-\log D)$, or equivalently $\Om^1(\log D)$, is a locally free $\O$-module, necessarily of rank $n$. 
\end{defi}

We will be concerned in this article with the following subclass of divisors.

\begin{defi}\label{2}
A free divisor $D$ is called \emph{linear} if $\Ga(\CC^n,\Der(-\log D))$ admits a basis $\delta_1,\dots,\delta_n$ such that each $\delta_i$ has linear coefficients with respect to the $\O$-basis $\p_1,\dots,\p_n$ of $\Der$ or equivalently each $\delta_i$ is homogeneous of degree zero with respect to the standard degree defined by $\deg x_i=1=-\deg\p_i$ on the variables and generators of $\Der$.
\end{defi}

Saito's criterion \cite[Thm.~1.8.(ii)]{Sai80} implies the following fundamental observation.

\begin{lemm}\label{28}
If $\delta_1,\dots,\delta_n$ is a basis of $\Ga(\CC^n,\Der(-\log D))$ for a linear free divisor $D$, then the homogeneous polynomial $\Delta=\det((\delta_i(x_j))_{i,j})\in\CC[x]$ of degree $n$ is a global defining equation for $D$.
\end{lemm}

Note that because $\Der(-\log D)$ can have no members of negative degree, $D$ cannot be isomorphic to the product of $\CC$ with a lower dimensional divisor. 
It turns out that linear free divisors are relatively abundant; the authors believe that in the current paper and in \cite{BM06}, recipes are given which allow the straightforward construction of more free divisors than have been described in the sum of all previous papers. 

\begin{exems}\label{1}\
\begin{asparaenum}
\item The normal crossing divisor $D=\{x_1\cdots x_n=0\}\subseteq\CC^n$ is a linear free divisor where
\[
x_1\p_1,\dots, x_n\p_n
\]
is a basis of $\Der(-\log D)$.
Up to isomorphism it is the only example among hyperplane arrangements, cf.~\cite[Ch.~4]{OT92}.

\item\label{1b} In the space $B_{2,3}$ of binary cubics, the discriminant $D$, which consists of binary cubics having a repeated root, is a linear free divisor. 
For $f(u,v)=xu^3+yu^2v+zuv^2+wv^3$ has a repeated root if and only if its Jacobian ideal does not contain any power of the maximal ideal $(u,v)$, and this in turn holds if and only if the four cubics
\[
u\p_uf,v\p_uf,u\p_vf,v\p_vf
\]
are linearly dependent. 
Writing the coefficients of these four cubics as the columns of the $4\times 4$ matrix
\[
A:=
\begin{pmatrix}
3x&0&y&0\\
2y&3x&2z&y\\
z&2y&3w&2z\\
0&z&0&3w
\end{pmatrix}
\]
we conclude that $D$ has equation $\det A=0$. 
After division by $3$ this determinant is
\[
-y^2z^2+4wy^3+4xz^3-18wxyz+27w^2x^2.
\]
In fact each of the columns of this matrix determines a vector field in $\Der(-\log D)$; 
for the group $\Gl_2(\CC)$ acts linearly on $B_{2,3}$ by composition on the right, and, up to a sign, the four columns here are the infinitesimal generators of this action corresponding to a basis of $\gl_2(\CC)$. 
Each is tangent to $D$, since the action preserves $D$. 

Further examples of irreducible linear free divisors can be found (though not under this name) in the paper \cite{SK77} of Sato and Kimura. 
Besides our example, two, of ambient dimension $12$ and $40$, are described in \cite[\S5, Prop.~11, 15]{SK77}, and by repeated application of castling transformations, cf.~\cite[\S2]{SK77}, it is possible to generate infinitely many more, of higher dimensions. 
\end{asparaenum}
\end{exems}

In Section~\ref{38} of this paper we describe a number of further examples of 
linear free divisors, and in Section~\ref{48} we prove some results about linear bases for the module $\Ga(\CC^n,\Der(-\log D))$, and go on to classify all linear free divisors in dimension $n\le 4$.

Linear free divisors provide a new insight into a conjecture of H.~Terao \cite[Conj.~3.1]{Ter78} relating the cohomology of the complement of certain divisors $D$ 
to the cohomology of the complex $\Om^\bu(\log D)$ of forms with logarithmic poles along $D$.
For linear free divisors, the link between the complex $\Ga(\CC^n,\Om^\bu(\log D))$ and $H^*(\CC^n\ssm D)$ can be understood as follows.

\begin{defi}\label{17}
For a linear free divisor $D$ defined by $\Delta\in\CC[x]$, we consider the subgroup
\[
G_D:=\{A\in\Gl_n(\CC)\mid A(D)=D\}=\{A\in\Gl_n(\CC)\mid\Delta\circ A\in\CC\cdot\Delta\}
\]
with identity component $G_D^\circ$ and Lie algebra $\gg_D$.
We call $D$ \emph{reductive} if $G_D^\circ$, or equivalently $\gg_D$, is reductive.
\end{defi}

It turns out that $\CC^n\ssm D$ is a single orbit of $G^\circ_D$ with finite isotropy group, so $H^*(\CC^n\ssm D;\CC)$ is isomorphic to the cohomology of $G^\circ_D$; this is explained in Section~\ref{8}. 
Moreover, $H^*(\Ga(\CC^n,\Om^\bu(\log D)))$ coincides with the Lie algebra cohomology of $\gg_D$ with complex coefficients.
For compact connected Lie groups $G$, a well-known argument shows that the Lie algebra cohomology coincides with the topological cohomology of the group. 
For linear free divisors the group $G^\circ_D$ is never compact, but the isomorphism also holds good for the larger class of reductive groups, and for a significant class of linear free divisors, $G^\circ_D$ is indeed reductive. 
In Section~\ref{15} we prove our main result:

\begin{theo}\label{3}
If $D$ is a reductive linear free divisor then
\beq\label{4}
H^*(\Ga(\CC^n,\Om^\bu(\log D)))\simeq H^*(\CC^n\smallsetminus D;\CC).
\eeq
\end{theo}

Among linear free divisors to which it applies are those arising as discriminants in representation spaces of quivers, as discussed in detail in \cite{BM06} and briefly in Section~\ref{31} below. 

Terao's conjecture remains open, though it has been answered in the affirmative for a very large class of arrangements in \cite{WY97}, using a technique developed in \cite{CMN96}. 
For general free divisors, a local result from which the global isomorphism of \eqref{4} follows holds when imposing the following additional hypothesis.

\begin{defi}\label{80}
A divisor $D$ is called \emph{quasihomogeneous} at $p\in D$ if the germ $(D,p)$ admits a local defining equation $\Delta\in\O_p$ that is weighted homogeneous with respect to weights $w_1,\dots,w_n\in\QQ_+$ in some local coordinate system $x_1,\dots,x_n$ centred at $p$.
Dividing $w_1,\dots,w_n$ by the weighted degree of $\Delta$, note that the preceding condition means that $\chi(\Delta)=\Delta$ where $\chi=\sum_{i=1}^nw_ix_i\p_i\in\Der(-\log D)_p$.
$D$ is called \emph{locally quasihomogeneous} if it is quasihomogeneous at $p$ for all $p\in D$.
We say \emph{homogeneous} instead of \emph{quasihomogeneous} if $w=1,\dots,1$.
\end{defi}

\begin{theo}[\cite{CMN96}]\label{6}
Let $D\subseteq\CC^n$ be a locally quasihomogeneous free divisor, let $U=\CC^n\smallsetminus D$, and let $j:U\to\CC^n$ be inclusion.
Then the de~Rham morphism
\beq\label{7}
\Om_X^{\bullet}(\log D)\to\mathbf{R}j_*\CC_U
\eeq
is a quasi-isomorphism.
\end{theo}

Grothendieck's \emph{Comparison Theorem} \cite{Gro66} asserts that a similar quasi-isomorphism holds for any divisor $D$, if instead of logarithmic poles we allow meromorphic poles of arbitrary order along $D$. 
Because of this similarity, we refer to the quasi-isomorphism of \eqref{7} as the Logarithmic Comparison Theorem (LCT) and to the global isomorphism \eqref{4} as the Global Logarithmic Comparison Theorem (GLCT).
Several authors have further investigated the range of validity of LCT, and established interesting links with the theory of $\D$-modules, in particular in
\cite{CN02}, \cite{CU05}, \cite{GS06}, \cite{Tor04}, and \cite{Wal05}.

Local quasihomogeneity was introduced in \cite{CMN96} as a technical device to make possible an inductive proof of the isomorphism in~\ref{6}.
Subsequently it turned out to have a deeper connection with the theorem. 
In particular by \cite{CCMN02}, for plane curves the logarithmic comparison theorem holds if and only if all singularities are quasihomogeneous.
The situation in higher dimensions remains unclear. 
There is as yet no counterexample to the conjecture that LCT is equivalent to the following weaker condition.

\begin{defi}\label{81}
A divisor $D$ is called \emph{Euler homogeneous} at $p\in D$ if there is a germ of vector field $\chi\in\mm_p\cdot\Der_p$ such that $\chi(\Delta)=\Delta$ for some local defining equation $\Delta\in\O_p$ of the germ $(D,p)$.
In this case, $\chi$ is called an \emph{Euler vector field} for $D$ at $p$.
$D$ is called \emph{strongly Euler homogeneous} if it is Euler homogeneous at $p$ for all $p\in D$.
\end{defi}

\begin{rema}
The Euler homogeneity of $D$ is independent of the choice of an equation. 
If $\chi$ is an Euler vector field at $p$ for $D$ defined by $\Delta\in\O_p$, and $u\in\O_p^*$ is a unit, then the defining equation $u\Delta$ of $D$ at $p$ satisfies an equation
\[
(\chi(u)+u)^{-1}u\chi(u\Delta)=(\chi(u)+u)^{-1}(\chi(u)+u)u\Delta=u\Delta
\]
with Euler vector field $(\chi(u)+u)^{-1}u\chi$.
\end{rema}

In Section \ref{5} we examine the examples described in Sections \ref{38} and \ref{48} with respect to local quasihomogeneity and strong Euler homogeneity. 
It turns out that all linear free divisors in dimension $n\le4$ are locally quasihomogeneous and there is no linear free divisor which we know not to be strongly Euler homogeneous.
The optimistic reader could therefore conjecture that all linear free divisors are strongly Euler homogeneous, and also fulfil LCT and so also GLCT.
We do not know any counter-example to these statements.

In Subsection~\ref{62} we give examples of quivers $Q$ and dimension vectors $\d$ for which the discriminant in $\Rep(Q,\d)$ is a linear free divisor but is not locally quasihomogeneous. 
In such cases Theorem~\ref{6} therefore does not apply, but Theorem~\ref{3} does.

In Subsection~\ref{65}, we show that a linear free divisor does not need to be reductive for LCT to hold.
However we do not know whether reductiveness of the group implies LCT for linear free divisors.
The property of being a linear free divisor is not local, and our proof of GLCT here is quite different from the proof of LCT in \cite{CMN96}. 

The fact that linear free divisors in $\CC^n$ arise as the complement of the open orbit of an $n$-dimensional connected algebraic subgroup of $\Gl_n(\CC)$, means that there is some overlap between the topic of this paper and of the paper \cite{SK77}, where Sato and Kimura classify irreducible \emph{prehomogeneous vector spaces}, that is, triples $(G,\rho,V)$, where $\rho$ is an irreducible representation of the algebraic group $G$ on $V$, in which there is an open orbit. 
However, the hypothesis of irreducibility means that the overlap is slight. 
Any linear free divisor arising as the complement of the open orbit in an irreducible prehomogeneous vector space is necessarily irreducible by \cite[\S4, Prop.~12]{SK77}, whereas among our examples and in our low-dimensional classification (in Section~\ref{48}) all the linear free divisors except one (Example~\ref{1}(\ref{1b})) are reducible.
Even where $G$ is reductive, the passage from irreducible to reducible representations in this context is by no means trivial, including as it does substantial parts of the theory of representations of quivers.

\section{Linear free divisors and subgroups of $\Gl_n(\CC)$}\label{8}

A degree zero vector field $\delta\in\Der$ can be identified with an $n\times n$ matrix $A=(a_{i,j})_{i,j}\in\CC^{n\times n}$ by $\delta=\sum_{i,j}x_ia_{i,j}\p_j=xA\p^t$.
Under this identification, the commutator of square matrices corresponds to the Lie bracket of vector fields.

Let $D\subseteq\CC^n$ be a reduced divisor defined by a homogeneous polynomial $\Delta\in\CC[x]$ of degree $d$.

\begin{defi}\label{23}
We denote by 
\[
\L_D:=\{xA\p^t\mid xA\p^t(\Delta)\in\CC\cdot\Delta\}\subseteq\Ga(\CC^n,\Der(-\log D))
\]
the Lie algebra of degree zero global logarithmic vector fields.
\end{defi}

Recall from Definition~\ref{2}, that $D$ is linear free if $\L_D$ contains a basis of $\Der(-\log D)$, and recall $G_D^\circ$ from Definition~\ref{17}.

\begin{lemm}\label{9}
$G_D^\circ$ is an algebraic subgroup of $\Gl_n(\CC)$ and $\gg_D=\{A\mid xA^t\p^t\in\L_D\}$.
\end{lemm}

\begin{proof} 
Clearly $G_D$ is a subgroup of $\Gl_n(\CC)$ and defined by a system of polynomial (determinantal) equations.
Thus $G_D$ and hence also $G_D^\circ$ is an algebraic subgroup of $\Gl_n(\CC)$.
The Lie algebra of $G_D^\circ$ consists of all $n\times n$-matrices $A$ such that
\[
\Delta\circ(I+A\eps)=a(\eps)\cdot\Delta\in\CC[\eps]\cdot\Delta
\]
where $\CC[\eps]=\CC[t]/\ideal{t^2}\ni[t]=:\eps$.
Taylor expansion of this equation with respect to $\eps$ yields
\[
\Delta+\p(\Delta)\cdot A\cdot x^t\cdot\eps=(a(0)+a'(0)\cdot\eps)\cdot\Delta
\]
and hence $a(0)=1$ and, by transposing the $\eps$-coefficient, $xA^t\p^t\in\L_D$.
The argument can be reversed to prove the converse by setting
\[
a(\eps):=1+(xA^t\p^t(\Delta)/\Delta)\cdot\eps.
\]
\end{proof}

\begin{lemm}\label{10}
The complement $\CC^n\ssm D$ of a linear free divisor is an orbit of $G_D^\circ$ with finite isotropy groups.
\end{lemm}

\begin{proof}
For $p\in\CC^n$, the orbit $G_D^\circ\cdot p$ is a smooth locally closed subset of $\CC^n$ whose boundary is a union of strictly lower dimensional orbits, cf.~\cite[Prop.~8.3]{Hum75}.
The orbit map $G_D^\circ\to G_D^\circ\cdot p$ sends $I_n+A\eps$ to $p+pA^t\eps$ and induces a tangent map
\beq\label{11}
\gg_D\onto T_p(G_D^\circ\cdot p),\quad A\mapsto pA^t.
\eeq
For $p\not\in D$, $\Der(-\log D)(p)$ and hence also $\L_D(p)$ is $n$-dimensional.
Then by Lemma \ref{9} and \eqref{11} $T_pG_D^\circ\cdot p$ and hence $G_D^\circ\cdot p$ are $n$-dimensional which implies the finiteness of the isotropy group of $p$ in $G_D^\circ$.
As this holds for all $p\not\in D$, the boundary of $G_D^\circ\cdot p$ must be $D$ and then $G_D^\circ\cdot p=\CC^n\ssm D$.
\end{proof}

Reversing our point of view we might try to find algebraic subgroups $G\subseteq\Gl_n(\CC)$ that define linear free divisors. 
This requires by definition that $G$ is $n$-dimensional and connected and by Lemma \ref{10} that there is an open orbit.
The complement $D$ is then a candidate for a free divisor.
Indeed $D$ is a divisor: comparing with \eqref{11}, $D$ is defined by the \emph{discriminant determinant}
\[
\Delta=\det\begin{pmatrix}A_1x^t&\cdots&A_nx^t\end{pmatrix}
\]
where $A_1,\dots,A_n$ is a basis of the Lie algebra $\gg$ of $G$ and we denote by $f=\Delta_\text{red}$ the reduced equation of $D$.
As the entries of the defining polynomial are linear, $\Delta$ is a homogeneous polynomial of degree $n$.
Thus, if $\Delta$ is not reduced, $D$ can not be linear free.
We shall see examples where this happens in the next section.
On the other hand, Saito's criterion \cite[Lem.~1.9]{Sai80} shows the following.

\begin{lemm}\label{12} 
Let the $n$-dimensional algebraic group $G$ act linearly on $\CC^n$ with an open orbit.
If $\Delta$ is reduced then $D$ is a linear free divisor.\qed
\end{lemm}

As a first step towards our main result, we now describe the cohomology of $\CC^n\ssm D$ in terms of $G_D^\circ$.

\begin{prop}\label{13}
Suppose that $D\subseteq\CC^n$ is a linear free divisor and let $G^\circ_{D,p}$ be the (finite) isotropy group of $p\in\CC^n\ssm D$ in $G_D^\circ$.
Then
\[
H^*(\CC^n\ssm D;\CC)=H^*(G^\circ_D;\CC)^{G^\circ_{D,p}}=H^*(G^\circ_D;\CC).
\] 
\end{prop}

\begin{proof}
By Lemma \ref{10}, $\CC^n\ssm D\cong G^\circ_D/G^\circ_{D,p}$ with finite $G^\circ_{D,p}$ and the first equality follows.
The second equality holds because $G^\circ_D$ is path connected, which means that left translation by $g\in G^\circ_{D,p}$ is homotopic to the identity and thus induces the identity map on cohomology.
\end{proof}

\begin{rema}\label{14} 
The argument for the second equality also shows that if $G^\circ_D$ is a finite quotient of the connected Lie group $G$ then $H^*(\CC^n\ssm D;\CC)\simeq H^*(G;\CC).$
We will use this below in calculating the cohomology of $\CC^n\ssm D$.
\end{rema}

\section{Cohomology of the complement and Lie algebra cohomology}\label{15}

Let $\gg$ be a Lie algebra. 
The complex of Lie algebra cochains with coefficients in the complex representation $V$ of $\gg$ has $k$th term $\bigwedge^k_{\CC}\Hom_{\CC}(\gg,V)\cong \Hom_{\CC}(\bigwedge^k_{\CC}\gg,V)$, and differential $d_L:\bigwedge^k_\CC\Hom(\gg,V)\to\bigwedge^{k+1}_\CC\Hom(\gg,V)$ defined by
\begin{gather}\label{16}
(d_L\om)(v_1\wedge\dots\wedge v_{k+1})=\\
\nonumber\sum_{i<j}(-1)^{i+j}\om([v_i,v_j]\wedge v_1\dots\wedge\widehat{v_i}\wedge\dots\wedge\widehat{v_j}\wedge\dots\wedge v_{k+1})+\\
\nonumber\sum_i(-1)^{i+1}v_i\cdot\om(v_1\wedge\dots\wedge\widehat{v_i}\wedge\dots\wedge v_{k+1}).
\end{gather}
The cohomology of this complex is the \emph{Lie algebra cohomology} of $\gg$ with coefficients in $V$ and will be denoted $H^*_A(\gg;V)$.

The exterior derivative of a differential $k$-form satisfies an 
identical formula:
\begin{gather}\label{71}
d\om(\chi_1\wedge\dots\wedge\chi_{k+1})=\\
\nonumber\sum_{i<j}(-1)^{i+j}\om([\chi_i,\chi_j]\wedge\chi_1\wedge\dots\wedge\widehat{\chi_i}\wedge\dots\wedge\widehat{\chi_j}\wedge\dots\wedge\chi_{k+1})+\\
\nonumber\sum_i(-1)^{i+1}\chi_i\cdot\om(\chi_1\wedge\dots\wedge\widehat{\chi_i}\wedge\dots\wedge\chi_{k+1}).
\end{gather}
Here the $\chi_i$ are vector fields. 

When $D$ is a free divisor and $V=\O_p$ for some $p\in D$, it is tempting to conclude from the comparison of \eqref{16} and \eqref{71} that the 
complex $\Om^\bu(\log D)$ coincides with the complex of Lie algebra cohomology, with coefficients in $\O_p$, of the Lie algebra $\Der(-\log D)_p$. 
For $\Om^1(\log D)_p$ is the dual of $\Der(-\log D)_p$, and $\Om^k(\log D)=\bigwedge^k\Om^1(\log D)$.
However, this identification is incorrect, since, in the complex $\Om^\bu(\log D)$, both exterior powers and $\Hom$ are taken over the ring of coefficients $\O$, rather than over $\CC$, as in the complex of Lie algebra cochains. 
The cohomology of $\Om^\bu(\log D)_p$ is instead the \emph{Lie algebroid} cohomology of $\Der(-\log D)_p$ with coefficients in $\O_p$.
Nevertheless, when $D$ is a linear free divisor, there is the following important
link between these two complexes.

Recall $\L_D$ from Definition~\ref{23}. 
\begin{lemm}\label{82}
Let $D$ be a linear free divisor.
The complex $\Ga(\CC^n,\Om^\bu(\log D))_0$ of global homogeneous differential forms of degree zero coincides with the complex $\bigwedge^\bu_\CC\Hom(\L_D,\CC)$ of Lie algebra cochains with coefficients in $\CC$. 
\end{lemm}

\begin{proof}
First we establish a natural isomorphism between the corresponding
terms of the two complexes. 
We have
\begin{align*}
\Om^1(\log D)&=\Hom_{\O}(\Der(-\log D),\O)\\
&=\Hom_{\O}(\L_D\otimes_{\CC}\O,\O)\\
&=\Hom_{\CC}(\L_D,\CC)\otimes_{\CC}\O.
\end{align*}
Since $\Hom_{\CC}(\L_D,\CC)$ is purely of degree zero, and the degree zero
part of $\O$ consists just of $\CC$, the degree zero part of $\Ga(\CC^n,\Om^1(\log D))$ is
\[
\Ga(\CC^n,\Om^1(\log D))_0=\Hom_{\CC}(\L_D,\CC).
\]
Since moreover $\Ga(\CC^n,\Om^1(\log D))$ has no part of negative degree, it follows that 
\[
\Ga(\CC^n,\Om^k(\log D))_0
=\Ga\bigl(\CC^n,\bigwedge^k_{\O}\Om^1(\log D)\bigr)_0 
=\bigwedge^k_{\CC}\Hom_{\CC}(\L_D,\CC).
\]
Next, we show that the coboundary operators are the same. 
Because we are working with constant coefficients, the second sum on the right in \eqref{16} vanishes. 
Let $\chi_1,\dots,\chi_{k+1}\in\L_D$.
Then for $\om\in\Ga(\CC^n,\Om^k(\log D))_0$ and $i\in\{1,\dots,k+1\}$, $\om(\chi_1\w\dots\w\widehat{\chi_i}\w\dots\w \chi_{k+1})$ is a constant. 
It follows that the second sum on the right in \eqref{71} vanishes. 
Thus, the coboundary operator $d_L$ and the exterior derivative $d$ coincide.
\end{proof}

More generally let us consider weights $w=w_1,\dots,w_n\in\QQ_+$ and assign the weight $w_i$ (resp.\ $-w_i$) to $x_i$ and $dx_i$ (resp.\ to $\p_i$). 
Then the set of homogeneous vector fields or differential forms of a given degree is well defined. 

\begin{lemm}\label{83}
Suppose that the divisor $D\subseteq\CC^n$ is quasihomogeneous with respect to weights $w=w_1,\dots,w_n\in\QQ_+$.   
Then the following holds for any open set $U\subseteq\CC^n$:
\begin{enumerate}
\item\label{83a} If $\om\in\Ga(U,\Om^k(\log D))$ is $w$-homogeneous, then $L_\chi(\om)=\deg_w(\om)\om$, where $L_\chi$ is the Lie derivative with respect the Euler vector field $\chi=\sum_{i=1}^nw_ix_i\p_i$.
\item\label{83b} For any closed $\om\in \Ga(U,\Om^k(\log D))$ with decomposition $\om=\sum_{j\geq j_0}\om_j$ into $w$-homogeneous parts, $\om-\om_0$ is exact.
\item\label{83c} If $\Ga(U,\Om^k(\log D))_r\subseteq\Ga(U,\Om^k(\log D))$ denotes the subspace of $w$-homogeneous forms of $w$-degree $r$, then
\[
\Ga(U,\Om^\bu(\log D))_0\into \Ga(U,\Om^\bu(\log D))
\]
is a quasi-isomorphism.
\end{enumerate}
\end{lemm}

\begin{proof}\
\begin{asparaenum}\pushQED{\qed}
\item is a straightforward calculation, using Cartan's formula
$L_{\chi}(\om)=d\iota_{\chi}\om+\iota_{\chi}d\om$, where $\iota_{\chi}$ is
contraction by $\chi$. 
\item follows, for, if $\om$ is closed, so is $\om_j$ for
every $j$, and thus
\[
\om-\om_0=
\sum_{0\neq j\geq j_0}\om_j=
L_\chi\bigl(\sum_{0\neq j\geq j_0}\frac{\om_j}{j}\bigr)=
d(\iota_\chi\bigl(\sum_{0\neq j\geq j_0}\frac{\om_j}{j}\bigr)).
\]
\item is now an immediate consequence.\qedhere
\end{asparaenum}
\end{proof}

From Lemma~\ref{82} and Lemma~\ref{83}\eqref{83c} applied to $U=\CC^n$ we deduce the following

\begin{prop}\label{84}
Let $D\subseteq\CC^n$ be a linear free divisor. 
Then 
\[\pushQED{\qed}
H^*(\Ga(\CC^n, \Om^\bu(\log D)))\cong H^*_A(\L_D;\CC).\qedhere
\]
\end{prop}

Recall $G_D^\circ$ and $\gg_D$ from Definition~\ref{17}.
From Propositions~\ref{13} and \ref{84} we deduce

\begin{coro}\label{18}
The global logarithmic comparison theorem holds for a linear free divisor $D$ if and only if
\beq\label{19}\pushQED{\qed}
H^*(G_D^\circ;\CC)\cong H^*_{A}(\gg_D;\CC).\qedhere
\eeq
\end{coro}

There is such an isomorphism if $G$ is a connected compact real Lie group with Lie algebra $\gg$ (which is not our situation here). 
Left translation around the group gives rise to an isomorphism of complexes
\[
T:\bigwedge^\bullet\gg^*\to\bigl(\Om^\bu(G)^G,d\bigr)
\]
where $\gg^*=\Hom_{\RR}(\gg,\RR)$ and $\Om^\bu(G)^G$ is the complex of left-invariant real-valued differential forms on $G$. 
Composing this with the inclusion 
\beq\label{20} 
\bigl(\Om^\bu(G)^G,d\bigr)\to \bigl(\Om^\bu(G),d\bigr)\eeq
and taking cohomology gives a morphism 
\beq\label{21}
\tau_G\colon H^*_A(\gg_D;\RR)\to H^*(G;\RR).
\eeq
If $G$ is compact, \eqref{21} is an isomorphism. 
For from each closed $k$-form $\om$ we obtain a left-invariant closed $k$-form $\om_A$ by averaging:
\[
\om_A:=\frac{1}{|G|}\int_G\ell_g^*(\om)d\mu_L,
\]
where $\mu_L$ is a left-invariant measure and $|G|$ is the volume of $G$ with respect to this measure.  
As $G$ is path-connected, for each $g\in G$ $\ell_g$ is homotopic to the identity, so $\om$ and $\ell_g^*(\om)$ are equal in cohomology. 
It follows from this that $\om$ and $\om_A$ are also equal in cohomology.

Of course, this does not apply directly in any of the cases discussed here, since $G_D$ is not compact. 
Nevertheless if $G_D$ is a reductive group, the complexified morphism \eqref{21} is an isomorphism.
We now briefly outline the necessary definitions. 
Let $G_0$ be a compact Lie group. 
Then (\cite[\S 5.4, Thm.~10]{OV90}) $G_0$ has a faithful real representation. 
It follows (\cite[\S 3.4, Thm.~5]{OV90}) that $G_0$ has an affine real algebraic group structure. 
This allows its complexification. 

\begin{defi}\label{25}\
\begin{enumerate} 
\item The complex Lie algebra representation is {\em reductive} if it is the direct sum of a semi-simple ideal and a diagonalizable ideal.
\item The complex linear algebraic group $G$ is {\em reductive} if its Lie algebra (representation) is reductive.
\end{enumerate}
\end{defi}

The term ``reductive'' is due to the fact that these groups are characterised, among complex algebraic groups, by the complete reducibility of every finite-dimensional complex representation.  
Chapter 5 of \cite{OV90} establishes a bijection between compact Lie groups and reductive complex linear algebraic groups: 

\begin{theo}[{\cite[\S5.2, Thm.~5]{OV90}}]
On any compact Lie group $K$ there exists a unique real algebraic group structure, whose complexification $K(\CC)$ is reductive.
Any reductive complex algebraic group possesses an algebraic compact real form (of which it is therefore the complexification).\qed
\end{theo} 

The significance of this notion for us derives from the following fact:

\begin{theo}[{\cite[\S5.2 Thm.~2]{OV90}}]\label{89}
Let $G$ be a complex reductive algebraic group with an $n$-dimensional compact real form $K$. 
Then $G$ is diffeomorphic to $K\times\RR^n$.\qed
\end{theo}

\begin{coro}\label{24} 
If $G$ is a connected reductive complex algebraic group with complex Lie algebra $\gg$ then 
\[
H^*_A(\gg;\CC)\simeq H^*(G;\CC).
\]
\end{coro}

\begin{proof}
Let $K$ be a compact real form of $G$.
By \ref{89}, inclusion of $K$ into $G=K(\CC)$ induces an ismorphism on cohomology. 
The Lie algebra $\gg$ of $G$ is the complexification of the Lie algebra $\mathfrak{k}$ of $K$, so we have 
\[
H^*_A(\gg;\CC)\simeq H^*(\mathfrak{k};\RR)\otimes\CC\simeq H^*(K;\RR)\otimes\CC\simeq H^*(K;\CC)\simeq H^*(G;\CC),
\]
where the second isomorphism comes from the isomorphism (\ref{21}).
\end{proof}

From Corollary~\ref{18}, Definition~\ref{25}, and Corollary~\ref{24} we now conclude Theorem ~\ref{3} as announced in the introduction: the Global Logarithmic Comparison Theorem holds for all reductive linear free divisors.

Using the reductiveness of the group $\Gl_n(\CC)$, we will show in the next section that the group $G_D$ is reductive for divisors obtained as discriminants in the representation spaces of quivers.
The subgroup $B_n\subseteq\Gl_n(\CC)$ of upper triangular matrices is not reductive, and appears as the group $G_D$ in Example~\ref{40} which shows that reductivity is not necessary for the GLCT to hold.

\section{Linear free divisors in quiver representation spaces}\label{31}

The following discussion summarises part of \cite{BM06}. 
A quiver $Q$ is a finite connected oriented graph; it consists of a set $Q_0$ of nodes and a set $Q_1$ of arrows joining some of them. 
For each arrow $\vp\in Q_1$ we denote by $t\vp$ (for ``tail'') and $h\vp$ (for ``head'') the nodes where it starts and finishes.
A (complex) representation $V$ of $Q$ is a choice of complex vector space $V_\alpha$ for each node $\alpha\in Q_0$ and linear map $V(\vp):V_{t\vp }\to V_{h\vp}$ for each arrow $\vp\in Q_1$.
For a fixed \emph{dimension vector}
\[
\d=(d_\alpha)_{\alpha\in Q_0}:=(\dim V_\alpha)_{\alpha\in Q_0}.
\]
and a choice of bases for the $V_\alpha$, $\alpha\in Q_0$, the \emph{representation space} of the quiver $Q$ of dimension $\d$ is
\[
\Rep(Q,\d):=\prod_{\vp\in Q_1}\Hom(\CC^{d_{t\vp}},\CC^{d_{h\vp}})\cong\prod_{\vp\in Q_1}\Hom(V_{t\vp},V_{h\vp}).
\] 
On this space the \emph{quiver group}
\[
\Gl(Q,\d):=\prod_{\alpha\in Q_0}\Gl_{d_\alpha}(\CC)\cong\prod_{\alpha\in Q_0}\Gl(V_{\alpha})
\]
acts, by
\beq\label{32}
\bigl((g_{\alpha})_{\alpha\in Q_0}\cdot V\bigr)_{\vp}:=g_{h\vp}V(\vp)g_{t\vp}^{-1}.
\eeq
This action factors through the group 
\beq\label{33}
Z:=\CC^*\cdot(I_{d_\alpha})_{\alpha\in Q_0}\subseteq\Z(\Gl(Q,\d))
\eeq
in the center of $\Gl(Q,\d)$ where $I_{d_\alpha}\in\Gl_{d_\alpha}(\CC)$ is the unit matrix.
The group $\Gl(Q,\d)/Z$ is reductive as, choosing a vertex $x_0\in Q_0$, we can consider it as a central quotient 
\beq\label{34}
\Gl(Q,\d)/Z\cong
\Bigl(\Sl_{d_{x_0}}(\CC)\times\prod _{x\in Q_0\ssm\{x_0\}}\Gl_{d_x}(\CC)\Bigr)\Big/\bigl(\mu_{d_{x_0}}\cdot \prod I_{d_x}\bigr)
\eeq
where $\mu_k\subseteq\CC^*$ denotes the cyclic subgroup of order $k$.
It acts faithfully on $\Rep(Q,\d)$.
For $\Rep(Q,\d)$ and $\Gl(Q,\d)/Z$ to play the role of $\CC^n$ and $G_D$ as in Section~\ref{8}, we must require
\begin{align}\label{35}
\sum_{n\in N}d_n^2-\sum_{\vp\in A}d_{t\vp}d_{h\vp}
&=\dim_{\CC}\Gl(Q,\d)-\dim_{\CC}\Rep(Q,\d)\\
\nonumber&=\dim Z
=1.
\end{align}
But this equality is not yet sufficient: 
it is also necessary that $\Gl(Q,\d)/Z$ has an open orbit.
This occurs if the general representation in $\Rep(Q,\d)$ is indecomposable. 
If both this last condition and \eqref{35} hold, $\d$ is called a \emph{real Schur root} of $Q$. 
In this case, there is a single open orbit, and the \emph{discriminant determinant} $\Delta$ defines its complement $D$, a divisor called the \emph{discriminant}. 
This is the consequence of a result due to Kraft and Riedtmann \cite[\S 2.6]{KR86}, which asserts that if the general representation is indecomposable it has only scalar endomorphisms.
Then
\beq\label{36}
\Gl(Q,\d)/Z\cong G_D=G_D^\circ.
\eeq
The above discussion combined with Theorem~\ref{3} proves the following 

\begin{theo}\label{37}
If $\d$ is a real Schur root of a quiver $Q$ and the discriminant $D$ in $\Rep(Q,\d)$ is reduced then $D$ is a linear free divisor that satisfies the GLCT.
\end{theo}

In \cite{BM06} it is shown that if, moreover, $Q$ is a \emph{Dynkin quiver}, i.e.\ its underlying unoriented graph is a Dynkin diagram of type $A_n$, $D_n$, $E_6$, $E_7$ or $E_8$, then $\Delta$ is always reduced, and thus defines a linear free divisor.
The significance of the Dynkin quivers is, that by a theorem of Gabriel \cite{Gab72}, they are the quivers of \emph{finite type}, i.e.\ the number of $\Gl(Q,\d)$ orbits in $\Rep(Q,\d)$ is finite. 
It is this that guarantees that $\Delta$ is always reduced, cf.~\cite[Prop.~5.4]{BM06}. 
It also implies that every root of a Dynkin quiver is a real Schur root.

\begin{coro}
If $\d$ is a (real Schur) root of a Dynkin quiver $Q$ then the discriminant $D$ in $\Rep(Q,\d)$ is a linear free divisor that satisfies GLCT.
\end{coro}

\begin{rema}
The argument showing that GLCT holds for the free divisors arising as discriminants in quiver representation spaces yields a simple topological proof of a theorem of V.~Kac \cite[p.~153]{Kac82} (see also \cite{Sch91}):
When $\d$ is a sincere (i.e.\ $\d_x>0$ for all $x\in Q_0$) real Schur root of a quiver $Q$ with no oriented cycles, the discriminant in $\Rep(Q,\d)$ has $|Q_0|-1$ irreducible components. 
The proof is this: 
the number of irreducible components of a divisor in a complex vector space is equal to the rank of $H^1$ of the complement.
From Theorem~\ref{3} we know that $H^1(\Rep(Q,\d)\ssm D;\CC)\simeq H^1_A(\gg_D;\CC)$; as by \eqref{34}
\[
\gg_D\simeq\sl_{d_{x_0}}(\CC)\oplus\bigoplus _{x\in Q_0\ssm\{x_0\}}\gl_{d_x}(\CC),
\]
it follows that
\[
H^1(\Rep(Q,\d)\ssm D;\CC)\simeq 0\oplus\bigoplus_{x\in Q_0\ssm\{x_0\}}H^1(\gl_{d_x}(\CC);\CC)
\]
and so has rank $|Q_0|-1$.

Another simple algebraic proof of Kac's theorem was pointed out to us by the referee.
It consists in determining the dimension of the vector space of rational function on $\CC^n$ with zeroes and poles along $D$ only and lifting them to the group $G_D$.
\end{rema}

\section{Examples of linear free divisors}\label{38}

The conclusion of Section~\ref{8} guides our search for linear free divisors. 
Our first example shows that the implication in Theorem~\ref{3} is not an equivalence.

We denote by
\beq\label{39}
E_{ij}=(\delta_{i,k}\cdot\delta_{j,l})_{k,l}\in\gl_n(\CC)
\eeq
the elementary matrix with $1$ in the $i$th row and $j$th column and $0$ elsewhere.

\subsection{A non-reductive example satisfying GLCT}

\begin{exem}\label{40}
For $n\ge2$, the group $B_n$ of $n\times n$ invertible upper triangular matrices is not reductive.
It acts on the space $\Sym_n(\CC)$ of symmetric matrices by transpose conjugation:
\[
B\cdot S=B^tSB.
\]
Under the corresponding infinitesimal action, the matrix $b$ in the Lie algebra $\bb_n$ gives rise to the vector field $\chi_b$ defined by
\[
\chi_b(S)=b^tS+Sb.
\]
The dimensions of $B_n$ and $\Sym_n(\CC)$ are equal. 
The discriminant determinant $\Delta$ is reduced and defines a linear free divisor $D=V(\Delta)$. 
\end{exem}

To see this, consider an elementary matrix $E_{ij}\in\bb_n$ and let $\chi_{ij}$ be the corresponding vector field on $\Sym_n(\CC)$. 
If $I$ is the $n\times n$ identity matrix, then $\chi_{ij}(I)=E_{ji}+E_{ij}$. 
The vectors $\chi_{ij}(I)$ for $1\leq i\leq j\leq n$ are therefore linearly independent, and $\Delta(I)\neq 0$.

For an $n\times n$ matrix $A$, let $A_j$ be the $j\times j$ matrix obtained by deleting the last $n-j$ rows and columns of $A$, and let $\det_j(A)=\det(A_j)$. 
If $B\in B_n$ and $S\in\Sym_n(\CC)$, then because $B$ is upper triangular, $(B^tSB)_j=B^t_jS_jB_j,$ and so $\det_j(B^tSB)=\det_j(B_j)^2\det_j(S)$. 
It follows that the hypersurface $D_j:=\{\det_j=0\}$ is invariant under the action, and the infinitesimal action of $B_n$ on $\Sym_n(\CC)$ is tangent to each. 
Thus $\Delta$ vanishes on each of them.
The sum of the degrees of the $D_j$ as $j$ ranges from $1$ to $n$ is equal to $\dim\Sym_n(\CC)$, and so coincides with the degree of $\Delta$. 
Hence $\Delta$ is reduced, and we conclude, by Lemma~\ref{12}, that 
$D=D_1\cup\cdots\cup D_n$ is a linear 
free divisor. 
In particular, when $n=2$, $D\subseteq\Sym_2(\CC)=\CC^3$ is the union of a quadric cone and one of its tangent planes.

We now give a proof that GLCT holds for $D$, in the spirit of the proofs of the preceding section, even though $D$ is not reductive. 
In fact LCT already follows, by Theorem~\ref{6}, from local quasihomogeneity, which we prove in Subsection~\ref{65} below.

\begin{prop}\label{73}
GLCT holds for the discriminant $D$ of the action of $B_n$ on $\Sym_n(\CC)$ in Example \ref{40}.
\end{prop}

\begin{proof}
The group $G_D^\circ$ is a finite quotient of the group $B_n$ of upper-triangular matrices in $\Gl_n(\CC)$. 
There is a deformation retraction of $B_n$ to the maximal torus $T$ consisting of its diagonal matrices, and, with respect to the standard coordinates $a_{ij}$ on matrix space, it follows that $H^*(B_n)$ is isomorphic to the free exterior algebra on the forms $da_{ii}/a_{ii}$. 
Each of these is left-invariant, and it follows that the map $\tau_{B_n}\colon H^*_A(\bb_n;\CC)\to H^*(B_n;\CC)$ from \eqref{21} is an epimorphism. 

Similarly, the Lie algebra complex $\bigwedge^\bu\bb_n^*$ has a contracting homotopy to its semisimple part. 
We may consider it as the complex of left-invariant forms on the group $B_n$. 
Assign weights $w_1,\ldots,w_n$ to the columns and weights $-w_1,\ldots,-w_n$ to the rows. 
This gives the elementary matrix $E_{ij}\in\bb_n$ the weight $w_i-w_j$.
If $\eps_{i,j}\in\bb_n^*$ denotes the dual basis and we assign the weight $0$ to $\CC$ then $\wt(\eps_{i,j})=-\wt(E_{i,j})$.
With respect to the resulting gradings of $\bb_n$ and $\bb_n^*$, both the Lie bracket and the differential $d_L$ of the complex $\bigwedge^\bu\bb_n^*$ are homogeneous of degree 0, cf.~\eqref{16}.

Let $E=\sum_iw_iE_{ii}$, and let $\iota_E:\bigwedge^\bu\bb_n^*\to\bigwedge^\bu\bb_n^*$ be the operation of \emph{contraction by $E$} defined by
\[
(\iota_E\om)(v_1\wedge\cdots\wedge v_k):=\om(E\wedge v_1\wedge\cdots\wedge v_k).
\]
Observe that for each generator $E_{ij}\in\bb_n$ we have
\beq\label{78}
[E,E_{ij}]=(w_i-w_j)\cdot E_{ij}=\wt(E_{ij})\cdot E_{ij}.
\eeq
We claim that the operation
\[
L_E:=\iota_Ed_L+d_L\iota_E,
\]
of taking the Lie derivative along $E$ has the effect of multiplying each homogeneous element of $\bigwedge^\bu\bb_n^*$ by its $w$-degree. 
Indeed the operation $L_E$ is a derivation of degree zero on $\bigwedge^\bu\bb_n^*$, 
and the result on 1 forms,
\[
L_E(\eps_{i,j})=(w_j-w_i)\eps_{i,j},
\]
is therefore sufficient and can be easily checked by direct calculation.

Thus $L_E$ defines a contracting homotopy from $\bigwedge^\bu\mathfrak{b}_n^*$ to its $w$-degree $0$ part $\bigl(\bigwedge^\bu \mathfrak{b}_n^*\bigr)_0$, by exactly the same calculation as in Lemma \ref{83}, but with $\Ga(U,\Om^k(\log D))$ and $L_\chi$ replaced respectively by $\bigwedge^\bu\mathfrak{b}_n^*$ and $L_E$. 
If we choose $w_1<\cdots<w_n$ then all off-diagonal members of the basis $\{\eps_{i,j}\}_{1\le i\le j\leq n}$ of $\bb_n^*$ have strictly positive $w$-degree.
It follows that 
\[
\bigwedge^\bu\mathfrak{b}_n^*\simeq\Bigl(\bigwedge^\bu\mathfrak{b}_n^*\Bigr)_0=\bigwedge^\bu\ideal{\eps_{1,1},\dots,\eps_{n,n}}=\bigwedge^\bu\ttt^*
\]
where $\ttt$ is the Lie algebra of the torus $T$ above.
The differential $d_L$ is zero on this subcomplex, showing that $\tau_{B_n}$ is an isomorphism.
\end{proof}

\subsection{Discriminants of quiver representations}

The following example, due to Ragnar-Olaf Buchweitz, is of the type discussed in Section~\ref{31}.

\begin{exem}\label{41}
In the space $M_{n,n+1}(\CC)$ of $n\times(n+1)$ matrices, let $D$ be the divisor defined by the vanishing of the product of
the maximal minors. 
That is, for each matrix $A\in M_{n,n+1}(\CC)$, let $A_j$ be $A$ minus its $j$'th column, and let $\Delta_j(A)=\det(A_j)$. 
Then 
\[
D=\{A\in M_{n,n+1}(\CC):\delta=\prod_{i=1}^{n+1}{\Delta}_j(A)=0\}.
\]
It is a linear free divisor. 
Here, as the group $G$ in Remark~\ref{14} we may take the product $Gl_n(\CC)\times \{\diag(1,\lambda_1,\dots,\lambda_n):\lambda_1,\dots,\lambda_n\in\CC^*\}$, acting by
\[
\bigl(A,\diag(1,\lambda_1,\dots,\lambda_n)\bigr)\cdot M=
A\cdot M\cdot\diag(1,\lambda_1,\dots,\lambda_n)^{-1}
\]
The placing of the $1$ in the first entry of the diagonal matrices is
rather arbitrary; it could be placed instead in any other fixed position
on the diagonal.
That $D$ is a linear free divisor follows from the fact that
\begin{enumerate}
\item the complement of $D$ is a single orbit, so the discriminant determinant
is not identically zero, and
\item the degree of $D$ is equal to the dimension of $G_D$, so the discriminant
determinant is reduced.
\end{enumerate}
\end{exem}

\begin{figure}[h]
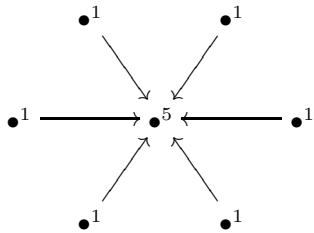

\caption{A star quiver with $1$ sink and $6$ sources and $\d=(5,1,\dots,1)$}\label{29}
\[
\xymat@C=.15in{&\bullet^1\ar[dr]&&\bullet^1\ar[dl]\\
\bullet^1\ar[rr]&&\bullet^5&&\bullet^1\ar[ll]\\
&\bullet^1\ar[ur]&&\bullet^1\ar[ul]}
\]
\end{figure}

In our Example~\ref{41}, $M_{n,n+1}(\CC)$ is the representation space of the star quiver consisting of one sink and $n+1$ sources, with dimension vector assigning dimension $n$ to the sink and $1$ to each of the sources.
The case $n=5$ is shown in Figure~\ref{29}.
Once we have chosen a basis for each $V_\alpha$, each $V(\vp)$ is represented by an $n\times 1$ matrix; together they make up an $n\times (n+1)$ matrix. 
So the basis identifies $\Rep(Q,\d)=M_{n,n+1}(\CC)$ and then
\beq\label{42}
\Gl(Q,\d)=\Gl_n(\CC)\times\Gl_1(\CC)^{n+1}
\eeq
and the action in \eqref{32} becomes
\beq\label{43}
(A,\lambda_1,\dots,\lambda_{n+1})\cdot M=AM\diag(\lambda_1^{-1},\dots,\lambda_{n+1}^{-1}).
\eeq
From \eqref{33}, \eqref{34}, and \eqref{36}, result isomorphisms
\beq\label{70}
\Gl(Q,\d)\cong G_D\times Z,\quad Z=\CC^*\cdot(I_n,(1,\dots,1)),
\eeq
defined by normalizing an arbitrary element in the second factor.

Many more examples of linear free divisors can be found by similar means in representation spaces of quivers. 
The next example, also from \cite{BM06}, is a non Dynkin quiver where finiteness of orbits fails and $\Delta$ is not reduced.

\begin{exem}\label{44}
Consider the star quiver of Example~\ref{41} with $n=3$ with $\d=(3,1,1,1,1)$, as before, and now reverse the direction of one of the arrows.
The result is shown in Figure~\ref{50}.
\begin{figure}[h]
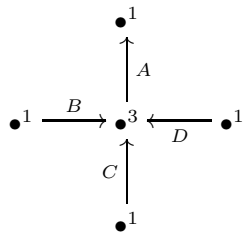

\caption{A modified star quiver with $\d=(3,1,1,1,1)$}\label{50}
\[
\xymat{&\bullet^1\\
\bullet^1\ar[r]^B&\bullet^3\ar[u]_A&\bullet^1\ar[l]^D\\
&\bullet^1\ar[u]^C}
\]
\end{figure}
The four hypersurfaces $\det(AB)=0$, $\det(AC)=0$, $\det(AD)=0$, $\det(BCD)=0$, are invariant under the action of the subgroup $G_D\subseteq\Gl(Q,\d)$ of Example~\ref{41}. 
However, the last of these is made up of infinitely many orbits: 
if the images of $B$, $C$ and $D$ lie in a plane $P$, then together with $\ker(A\cap P)$ they determine a cross-ratio. 
The discriminant determinant is equal, up to a scalar factor, to 
\[
\Delta=\det(AB)\cdot\det(AC)\cdot\det(AD)\cdot(\det(BCD))^2.
\]
\end{exem}

\subsection{Incomplete collections of maximal minors}\label{45} 

In the space $M_{m,n}(\CC)$ of $m\times n$ matrices with $n>m+1$, the product of all of the maximal minors no longer defines a linear free divisor, by reason of its degree. 
However, certain collections of $n$ maximal minors do define free divisors. There is a simple procedure for generating infinitely many  such collections, first described in \cite{Nie05}:

The space $M_{m,n}(\CC)$ can still be viewed as $\Rep(Q,\d)$ where $Q$ is the star quiver of Example~\ref{41} with $1$ sink and $n$ sources, and $\d=(m,1,\dots,1)$. 
As before, the quiver group $\Gl(Q,\d)$ acts with $1$-dimensional kernel $Z$, but now
\[
\dim\Gl(Q,\d)-1=m^2+n-1<mn=\dim M_{m,n}(\CC),
\]
making an open orbit impossible.
Therefore we replace $\Gl(Q,\d)$ by a group
\beq\label{79}
G:=\Gl_m(\CC)\times G_R
\eeq
with $\dim G/Z=\dim M_{m,n}$ by augmenting the second factor in \eqref{42} to a group $G_R\subseteq\Gl_n(\CC)$ with $\dim  G_R=mn-m^2+1$. 
To construct $G_R$, we consider an auxiliary quiver $\tilde Q=(Q_0,Q_1)$ with $Q_0=\{1,\dots,n\}$ and $Q_1\subseteq Q_0^2$ satisfying the following conditions:
\begin{itemize}
\item $|Q_1|=mn-m^2+1$;
\item $(i,i)\in Q_1$ for all $i\in Q_0$;
\item $(i,j)\in Q_1$ and $(j,k)\in Q_1$ implies that $(i,k)\in Q_1$.
\end{itemize}
These conditions are exactly those we need for the following formula:
\beq\label{72}
G_R:=\CC^{Q_1}\ssm\{\det=0\}\subseteq\CC^{n\times n}\ssm\{\det=0\}=\Gl_n(\CC)
\eeq
to define a group.
We write $(Q_0,Q_1)=:Q(G_R)$.
This group $G_R$ is generated by $\Gl_1(\CC)^n=\diag(\CC^*,\dots,\CC^*)$ and $mn-m^2-n+1$ supplementary elementary matrices $I_n+\CC\cdot E_{i,j}$ with $(i,j)\in Q_1$ and $i\ne j$, cf.~\eqref{39}.

The action of $G$ on $M_{m,n}(\CC)$ extends that in \eqref{43} by right multiplication of $G_R$ and factors through $G/Z$ with $Z=\CC\cdot (I_m, I_n)$ which is, as in \eqref{34}, a central quotient
\beq\label{74}
G/Z\cong(\Sl_m(\CC)\times G_R)/\mu_m\cdot(I_m,I_n)
\eeq
where $\mu_k\subseteq\CC^*$ denotes the cyclic subgroup of order $k$. 

\begin{prop}\label{76}\
\begin{asparaenum}
\item If the discriminant determinant $\Delta$ of the action of $G$ is not identically zero and the action of $G$ preserves the divisors of zeros of precisely $n$ distinct $m\times m$ minors, then the union of these divisors is a linear free divisor $D=V(\Delta)$.
\item If the action of $G$ preserves the divisor of zeros of more than $n$ distinct $m\times m$ minors then $\Delta$ is identically zero.
\end{asparaenum}
\end{prop}

\begin{proof}
Any algebraic set preserved by the action of $G$ is contained in $V(\Delta)$.
By construction, if $\Delta$ is not identically zero then its degree is $mn$.
If moreover the action of $G$ preserves the zero set of $n$ distinct $m\times m$ minors then $\Delta$ is reduced and Lemma~\ref{12} shows that $V(\Delta)$ is a linear free divisor. 
\end{proof}

\begin{lemm}\label{46}
Right multiplication by $I_n+\CC\cdot E_{i,j}$ preserves the divisor defined by an 
$m\times m$ minor if and only if the minor either contains column $i$ of the generic matrix or does not contain column $j$.
\end{lemm}

\begin{proof}
Suppose that the $m\times m$ submatrix $M'$ of the generic $m\times n$ matrix $M$ contains column $j$ but not column $i$. 
Let $p$ be a point of $\{\det M'=0\}$ at which $M'$ has rank $m-1$ and the matrix $M''$ obtained from $M'$ by replacing column $j$ by column $i$ has rank $m$. 
Then $\det(M'\cdot(I_n+E_{i,j}))(p)\neq 0$. 
That is, $\cdot(I_n+\CC\cdot E_{i,j})$ does not preserve $\{\det M'=0\}$. 
Similarly, $\det M'$ is clearly invariant under $\cdot(I_n+\CC\cdot E_{i,j})$ if $M'$ contains both columns $i$ and $j$.
\end{proof}

\begin{exem}\label{47}
\begin{annalif}
Let $m=3$ and $n=6$. 
The quiver
\[
\xymatrix@C=3mm@R=2mm{
1\ar[dr]&&3\ar[dl]\ar[dr]&&5\ar[dl]\\
&2&&4&}
\]
determines admissible minors $M_{123}$, $M_{345}$, $M_{135}$, $M_{136}$, $M_{156}$, $M_{356}$, and the associated divisor (the zero locus of their product) is a linear free divisor. 
Other linear free divisors that can be constructed by these methods are listed, for small values of $m,n$, in the preprint version \cite{GMNS06} of this article.
They are not in general reductive. 
\end{annalif}
\begin{arxiv}
\begin{asparaenum}

\item Take $m=2$, $n=4$. 
We refer to the $2\times 2$ submatrix of the generic matrix containing columns $i$ and $j$ as $M_{ij}$.
We need one supplementary column operation type. 
By reordering the columns we may assume this is $\cdot(I_n+\CC\cdot E_{3,4})$.
The set of submatrices satisfying the condition of Lemma~\ref{46} is $M_{12},M_{13},M_{23},M_{34}$. 
The product of their determinants therefore defines a linear free divisor. 
The group $G_R$ consists of matrices of the form 
\[
\begin{pmatrix}
1&0&0&0\\
0&\lambda_2&0&0\\
0&0&\lambda_3&\lambda_{34}\\
0&0&0&\lambda_4
\end{pmatrix}
\]

\item\label{47b}Take $m=2$, $n=5$. 
Then we require two supplementary generators which means that the quiver $Q(\tilde G_R)$ has two arrows.
Our Table~\ref{85} shows all such quivers with no more than five vertices (up to reordering of columns).
Note that the quiver
\[
\xymat{1\ar[r]&2\ar[r]&3}
\]
does not appear since the conditions on $Q_1$ force it to become 
\[
\xymat{1\ar[r]\ar@/^/[rr]&2\ar[r]&3}.
\]
The symbols $*$ and $\star$ denote arbitrary elements of $\CC^*$ and $\CC$ respectively. 

\begin{longtable}{|c|c|c|c|}
\caption{LFDs from collections of minors for $m=2$, $n=5$}\label{85}\\
\hline
\rule{0pt}{11pt}$Q(\tilde G_R)$ & $\tilde G_R$ & Admissible minors & LFD?\\
\hline
$
\xymat@R=1mm
{&4\\
3\ar[ur]\ar[dr]&\\
&5}
$
&
$
\begin{pmatrix}
*&0&0&0&0\\
0&*&0&0&0\\
0&0&*&\star&\star\\
0&0&0&*&0\\
0&0&0&0&*
\end{pmatrix}
$
&
$
M_{12},M_{13},M_{23},M_{34},M_{35}
$
&
Yes
\\
\hline
$\xymat@R=1mm
{2\ar[r]&3\\
4\ar[r]&5}
$
&
$
\begin{pmatrix}
*&0&0&0&0\\
0&*&\star&0&0\\
0&0&*&0&0\\
0&0&0&*&\star\\
0&0&0&0&*
\end{pmatrix}
$
&
$
M_{12},M_{14},M_{23},M_{24},M_{45}
$
&
Yes
\\
\hline
$
\xymat@R=1mm
{3\ar[dr]&\\
&5\\
4\ar[ur]&}
$
&
$
\begin{pmatrix}
*&0&0&0&0\\
0&*&0&0&0\\
0&0&*&0&\star\\
0&0&0&*&\star\\
0&0&0&0&*
\end{pmatrix}
$
&
$
M_{12},M_{13},M_{23},M_{14},M_{24},M_{34}
$
&
No
\\
\hline
$
\xymat@R=1mm
{1\ar@/^/[r]&\ar@/^/[l]2}
$
&
$
\begin{pmatrix}
*&\star&0&0&0\\
\star&*&0&0&0\\
0&0&*&0&0\\
0&0&0&*&0\\
0&0&0&0&*\\
\end{pmatrix}
$
&
$
M_{12},M_{34},M_{35},M_{45}
$
&
No
\\
\hline
\end{longtable}

Only the first and second yield linear free divisors. 
The fourth fails because although the group has the right dimension and an open orbit, the discriminant determinant is not reduced: 
the minor $M_{12}$ divides it with multiplicity $2$.

\item Take $m=3$, $n=5$. Again we need two supplementary generators, and there are the same four cases as when $(m,n)=(2,5)$. 
In each case the group $G_R$ is the same as in the previous example. 
The Table~\ref{86} of those which give rise to free divisors is different from the previous example now.

\begin{longtable}{|c|l|l|}
\caption{LFDs from collections of minors for $m=3$, $n=5$}\label{86}\\
\hline
\rule{0pt}{11pt}Admissible minors & LDF?\\
\hline
$M_{123},M_{134},M_{234},M_{345},M_{135},M_{235}$ & No\\
\hline
$M_{123},M_{124},M_{234},M_{145},M_{245}$ & Yes\\
\hline
$M_{123},M_{124},M_{134},M_{234},M_{345}$ & Yes\\
\hline
$M_{123},M_{124},M_{125},M_{345}$ & No\\
\hline
\end{longtable}

\item Take $m=3$, $n=6$. 
We need four supplementary generators and show all possible quivers with no more than six vertices in Table~\ref{87}. 

\begin{longtable}{|c|c|l|c|}
\caption{LFDs from collections of minors for $m=3$, $n=6$}\label{87}\\
\hline
\rule{0pt}{11pt}$Q(\tilde G_R)$ & Admissible minors & LFD?\\
\hline
\xymat@C=3mm@R=1mm 
{2&&1\ar[ll]\ar[dl]\ar[dr]\ar[rr]&&5\\
&3&&4&}
&
$
\begin{matrix}
M_{123},M_{124},M_{125},M_{134},M_{135},\\
M_{145},M_{126},M_{136},M_{146},M_{156}                         
\end{matrix}
$
&
No
\\
\hline
\xymat@C=3mm@R=1mm
{&2\ar[dr]&&3\ar[dl]&\\
1\ar[rr]&&5&&4\ar[ll]}
&
$
\begin{matrix}
M_{123},M_{124},M_{134},M_{234},M_{126},\\
M_{136},M_{146},M_{236},M_{246},M_{346}                         
\end{matrix}
$
&
No
\\
\hline
\xymat@C=3mm@R=1mm
{&1\ar[dl]\ar[dr]&&&4\ar[dl]\ar[dr]&\\
2&&3&5&&6}
&
$
\begin{matrix}
M_{123},M_{124},M_{134},M_{456},M_{145},M_{146}
\end{matrix}
$
&
Yes
\\
\hline
\xymat@C=3mm@R=1mm
{1\ar[dr]&&2\ar[dl]&4\ar[dr]&&5\ar[dl]\\
&3&&&6}
&
$
\begin{matrix}
M_{123},M_{124},M_{125},M_{145},M_{245},M_{456}
\end{matrix}
$
&
Yes
\\
\hline
\xymat@C=3mm@R=1mm
{1\ar[dr]&&2\ar[dl]&4&&5\\
&3&&&6\ar[ul]\ar[ur]}
&
$
\begin{matrix}
M_{123},M_{126},M_{146},M_{156},\\
M_{246},M_{256},M_{456}
\end{matrix}
$
&
No
\\
\hline
\xymat@C=6mm@R=2mm
{1\ar[d]&3\ar[dr]&4\ar[d]&5\ar[dl]\\
2&&6}
&
$
\begin{matrix}
M_{123},M_{124},M_{125},M_{134},\\
M_{135},M_{145},M_{345}
\end{matrix}
$
&
No
\\
\hline
\xymat@C=6mm@R=2mm
{1\ar[d]&&3\ar[dl]\ar[d]\ar[dr]\\
2&4&5&6}
&
$
\begin{matrix}
M_{123},M_{134},M_{135},M_{136},\\
M_{345},M_{346},M_{356}
\end{matrix}
$
&
No
\\
\hline
\xymat@C=3mm@R=2mm
{1\ar[dr]&&3\ar[dl]\ar[dr]&&5\ar[dl]\\
&2&&4&}
&
$
\begin{matrix}
M_{123},M_{345},M_{135},M_{136},M_{156},M_{356}
\end{matrix}
$
&
Yes
\\
\hline
\xymat@C=3mm@R=2mm
{&1\ar[dl]\ar[dr]&&4\ar[dl]\ar[dr]&\\
2&&3&&5}
&
$
\begin{matrix}
M_{124},M_{126},M_{134},M_{456},M_{145},M_{146}
\end{matrix}
$
&
Yes
\\
\hline
\xymat@R=3mm@C=6mm
{1\ar[d]&3\ar[d]\ar[dr]&4\ar[d]\\
2&5&6}
&
$
\begin{matrix}
M_{123},M_{124},M_{134},M_{135},M_{345},M_{346}
\end{matrix}
$
&
Yes
\\
\hline
\xymat@R=3mm@C=6mm
{1\ar[d]&3\ar[dl]\ar[d]\ar[dr]&\\
2&4&5}
&
$
\begin{matrix}
M_{123},M_{134},M_{135},M_{136},\\
M_{345},M_{346},M_{356}
\end{matrix}
$
&
No
\\
\hline
\xymat@R=3mm@C=6mm
{1\ar[d]\ar[dr]&2\ar[d]&3\ar[dl]\\
4&5}
&
$
\begin{matrix}
M_{123},M_{124},M_{134},M_{126},\\
M_{136},M_{146},M_{236}
\end{matrix}
$
&
No
\\
\hline
\xymat@R=3mm@C=6mm
{1\ar[d]\ar[r]&4\\
3&2\ar[l]\ar[u]}
&
$
\begin{matrix}
M_{123},M_{124},M_{125},M_{126},M_{156},M_{256}
\end{matrix}
$
&
Yes
\\
\hline
\xymat@R=3mm@C=6mm
{1\ar[d]\ar[r]&2\ar[dl]&\\
3&4\ar[r]&5}
&
$
\begin{matrix}
M_{123},M_{124},M_{126},M_{145},M_{146},M_{456}
\end{matrix}
$
&
Yes
\\
\hline
\xymat@R=3mm@C=6mm
{1\ar[d]\ar[r]&2\ar[dl]\\
3&4\ar[l]}
&
$
\begin{matrix}
M_{124},M_{125},M_{126},M_{145},\\
M_{146},M_{156},M_{456}
\end{matrix}
$
&
No
\\
\hline
\xymat@R=3mm@C=6mm
{&1\ar[dl]\ar[d]\ar[dr]&\\
2\ar[r]&3&4}
&
$
\begin{matrix}
M_{123},M_{124},M_{145},M_{146},\\
M_{125},M_{126},M_{156}
\end{matrix}
$
&
No
\\
\hline
\end{longtable}
\end{asparaenum}

\end{arxiv}
\end{exem}

\begin{prop}\label{75}
Let $D=V(\Delta)$ be a linear free divisor as constructed above.
If $Q(G_R)$ has no oriented loops then GLCT holds for $D$.
\end{prop}

\begin{proof}
As in the proof of Proposition~\ref{73}, one can show that
\[
\tau_{G_R}\colon H^*_A(\mathfrak{g}_R;\CC)\to H^*(G_R;\CC)
\]
from \eqref{21} is an isomorphism.
Here the absence of oriented loops serves as a replacement for the upper triangularity in the preceding proof.
Indeed, if there are no oriented loops in $Q( G_R)$, it is possible to order the
vertices of $Q(G_R)$, and thus the rows of the matrices in $M_{m,n}$, so that $i<j$ whenever there is an arrow from $i$ to $j$. 
This puts all of the matrices of $G_R$ into upper triangular form. 
It follows both that $G_R$ has a deformation retraction to its maximal torus $T$ consisting of diagonal matrices, and that the same contracting homotopy as in the proof of \ref{73} shows that the inclusion 
$\bigwedge^\bu\mathfrak{t}^*\to \bigwedge^\bu\mathfrak{g}_R^*$ is
a homotopy equivalence. 
Thus $H^*(T):H^*_A(\mathfrak{g}_R;\CC)\to H^*(G_R;\CC)$ is an isomorphism.

Also for $G=\Sl_m(\CC)$ the map $\tau_G$ from \eqref{21} is an isomorphism.
So by applying the the K\"unneth formulas for both Lie algebra and complex cohomology, the same holds for $G=\Sl_m(\CC)\times G_R$. 
By \eqref{74}, $G/Z$ is connected as a finite quotient of the connected group $\Sl_m(\CC)\times G_R$.
By Proposition~\ref{76} $G_D^\circ$ is then also a connected finite quotient of $G/Z$ hence of $\Sl_m(\CC)\times G_R$, and GLCT holds for $D$ by Corollary~\ref{18}.
\end{proof}

\section{Classification in small dimensions}\label{48}

\subsection{Structure of logarithmic vector fields}\label{49}

Let $\delta, \xi\in\Der$ and let $g\in\O$. 
To emphasise the action of $\delta$ on $\O$ and on $\Der$, in place of $dg(\delta)$ we write $\delta(g)$, and in place of $[\delta,\xi]$ we write $\delta(\xi)=\ad_\delta(\xi)$.
The degree $k$ parts of $\Ga(\CC^n,\Der)$ and $\Ga(\CC^n,\Der(-\log D))$ with respect to $\deg(x_i)=1=-\deg(\p_i)$ will be denoted by $\Ga(\CC^n,\Der)_k$ and $\Ga(\CC^n,\Der(-\log D))_k$ respectively.
For $\delta\in\Ga(\CC^n,\Der)_0$, we write $\delta_S$ for its semisimple part and $\delta_N$ for its nilpotent part.

Let $D\subseteq\CC^n$ be a linear free divisor defined by the homogeneous polynomial $\Delta=\det((\delta_i(x_j))_{i,j})\in\CC[x]$ of degree $n$ as in Lemma~\ref{28} where $\delta=\delta_1,\dots,\delta_n$ is a global degree $0$ basis of $\Der(-\log D)$.
Then $\delta_i(\Delta)\in\CC\cdot\Delta$ and there is the standard Euler vector field $\chi=\sum_ix_i\p_i\in\langle\delta_1,\dots,\delta_n\rangle_{\CC}$.
Since $\chi(\Delta)/\Delta=n\ne0$, we can assume that $\delta_1=\chi$ and $\delta_i(\Delta)=0$ for $i=2,\dots,n$.
So $\delta_2,\dots,\delta_n$ is a global degree $0$ basis of the annihilator $\Der(-\log\Delta)$ of $\Delta$ which is a direct factor of $\Der(-\log D)$.

Since $\chi$ vanishes only at the origin, the origin of the affine coordinate system $x=x_1,\dots,x_n$ is uniquely determined.
A coordinate change between two degree $0$ bases of $\Der(-\log D)$ can always be chosen linear. 
Among all possible linear coordinate changes, let $s+1$ be the maximal number of linearly independent diagonal logarithmic vector fields.

\begin{theo}\label{51}
There exists a global degree $0$ basis $\chi$, $\sigma_1,\dots,\sigma_s$, $\nu_1,\dots,\nu_{n-s-1}$ 
of $\Der(-\log D)$ such that
\begin{enumerate}
\item\label{51a} $\chi(\sigma_i)=0$ and $\chi(\nu_j)=0$,
\item\label{51b} the $\sigma_i$ are simultaneously diagonalizable with eigenvalues in $\QQ$ and $\sigma_i(\Delta)=0$,
\item\label{51c} the $\nu_j$ are nilpotent and $\nu_j(\Delta)=0$,
\item\label{51d} $\sigma_i(\nu_j)\in\QQ\cdot\nu_j$ and 
$\sum_j\sigma_i(\nu_j)/\nu_j+\tr(\sigma_i)=0$.
\item\label{51e} If $\delta\in\Ga(\CC^n,\Der(-\log D))_0$ with $\sigma_i(\delta)=0$ for $i=1,\dots,s$ then $\delta_S\in\ideal{\sigma_1,\dots,\sigma_s}_\CC$.
\end{enumerate}
Moreover, $s\ge 1$ and if $s=n-1$ then $\Delta=x_1\cdots x_n$ defines a normal crossing divisor.
\end{theo}

\begin{proof}
It is easy to check that the formal coordinate changes used in \cite{GS06} reduce to \emph{linear} coordinate changes in the case of \emph{linear} free divisors.
Thus \eqref{51a}-\eqref{51c}, \eqref{51e}, and the first part of \eqref{51d} follow from \cite[Thm.~5.4]{GS06}. 

For the second part of \eqref{51d}, we set $\delta_1,\dots,\delta_n=\chi,\sigma_1,\dots,\sigma_s,\nu_1\dots,\nu_{n-s-1}$ and rewrite $\Delta$ as
\beq\label{27}
\Delta=\sum_{\alpha\in S_n}\sgn(\alpha)\cdot\delta_1(x_{\alpha_1})\cdots\delta_n(x_{\alpha_n}).
\eeq
Let us choose coordinates in which all $\sigma _i$ are diagonal: $\sigma_i=\sum_j w_{i,j}x_j\p_j$. 
The equation $\sigma_i(\Delta)=0$ means that $\Delta$ is weighted homogeneous of degree zero when we assign to the variable $x_j$ the weight $w_{i,j}=\sigma_i(x_j)/x_j$.
The weighted degree of $\delta_j(x_k)$ is then $\sigma_i(\delta _j )/\delta_j+w_{i,k}$. 
This implies the second part of \eqref{51d}, since each term in the sum \eqref{27} has the same weighted degree $\sum_j\big(\sigma_i(\delta_j)/\delta_j+w_{i \alpha_j}\big)=\sum_j\sigma_i(\delta_j)/\delta_j+\tr(\sigma_i)$.

Now assume that $s=0$.
Then the vector space generated by the $\nu _i$ is entirely made of nilpotent elements and we can apply Engel's Theorem \cite[I.4]{Ser87}, and $\nu_1,\dots,\nu_{n-1}$ can be chosen upper triangular.
But then $\Delta$ is clearly divisible by the square of the first variable $x_1$ and hence not reduced.
So $s=0$ is impossible.

If $s=n-1$, then $\Delta$ must be a monomial and hence $\Delta=x_1\cdots x_n$ defines a normal crossing divisor.
\end{proof}

\begin{rema}
In Theorem~\ref{51}, one can perform the Gauss algorithm on the diagonals of $\sigma_1,\dots,\sigma_s$.
Then $\sigma_i\equiv x_i\p_i\mod\sum_{j=s+1}^n\CC\cdot x_j\p_j$.
\end{rema}

We shall frequently use the following simple fact.

\begin{lemm}\label{52}
Let $\sigma=\sum_{i}w_ix_i\p_i$.
Then $x_i\p_j$ is an eigenvector of $\ad_\sigma$ for the eigenvalue 
$w_i-w_j$.
\end{lemm}

\subsection{The case $s=n-2$}\label{53}

\begin{lemm}\label{54}
Let $s=n-2$ in the situation of Theorem~\ref{51}.
Then $-\sigma_k(\nu_1)/\nu_1=\tr(\sigma_k)\ne0$ for some $k$ and $\nu_1$ has, after normalization, two entries equal to $1$ and all other entries equal to $0$.
\end{lemm}

\begin{proof}
If $s=n-2$ then for any $\sigma \in \{\sigma_1,\dots,\sigma_{n-2}\}$, $\sigma(\nu_1)/\nu_1+\tr(\sigma)=0$.
Hence, a monomial $x_i\p_j$ in $\nu_1$ gives a relation $w_i-w_j+\sum_kw_k=0$ on the diagonal entries $w_1,\dots,w_n$ of $\sigma$.
Since $3$ of these relations with $i\ne j$ and also $\sigma_1,\dots,\sigma_{n-2}$ are linearly independent, $\nu_1$ can have at most $2$ nonzero nondiagonal entries.
If $\sigma_k(\nu_1)/\nu_1\ne0$ for some $k$ then $\nu_1$ is strictly triangular with at most $2$ nonzero entries after ordering the diagonal of $\sigma_k$.
If $\nu_1$ has only one nonzero entry then $\Delta$ is divisible by the square of a variable, a contradiction.
Both nonzero entries of $\nu_1$ can be normalized to $1$.
If $\sigma_k(\nu_1)/\nu_1=0$ for all $k$ then the nonzero entries of $\nu_1$ are in a $2$-dimensional simultaneous eigenspace of $\chi,\sigma_1,\dots,\sigma_{n-2}$.
Otherwise, there are $3$ linearly independent relations $w_{i_1}=w_{j_1}$, $w_{i_2}=w_{j_2}$, $\sum_kw_k=0$ on the diagonal entries of $\sigma_1,\dots,\sigma_{n-2}$, a contradiction to the linear independence of these vector fields.
But then $\nu_1$ has only one nonzero entry after a linear coordinate change, a contradiction as before.
\end{proof}

To simplify the notation, we shall write $\equiv$ for equivalence modulo $\CC^*$.
By Lemma~\ref{54}, we may assume that $\nu_1=x_k\p_1+x_l\p_2$ where $1\ne k\ne l\ne 2$.
Then 
\[
\Delta=
\begin{vmatrix}
x_1 & x_2 & x_3 & \cdots & x_n \\
a_{2,1}x_1 & a_{2,2}x_2 & a_{2,3}x_3 & \cdots & a_{2,n}x_n \\
\vdots & \vdots & \vdots && \vdots \\
a_{n-1,1}x_1 & a_{n-1,2}x_2 & a_{n-1,3}x_3 & \cdots & a_{1,n}x_n \\
x_k & x_l & 0 & \cdots & 0
\end{vmatrix}
\equiv(x_2x_k-x_1x_l)x_3\cdots x_n.
\]

As $\Delta$ is reduced, there are, up to coordinate changes, only two non-normal-crossing cases:

\subsubsection{$k=2$, $l=3$}\label{55}

Then $\Delta$ comes from the linear free divisor of Example~\ref{40} in dimension $3$:
\[
\Delta=(x_2^2-x_1x_3)x_3\cdots x_n\equiv
\begin{vmatrix}
x_1 & x_2 & x_3 \\
4x_1 & x_2 & -2x_3 \\
2x_2 & x_3 & 0
\end{vmatrix}
\cdot x_4\cdots x_n.
\]

\subsubsection{$k=3$, $l=4$}

Then $\Delta$ comes from a linear free divisor in dimension $4$:
\[
\Delta=(x_2x_3-x_1x_4)x_3\cdots x_n\equiv
\begin{vmatrix}
x_1 & x_2 & x_3 & x_4 \\
x_1 & 2x_2 & -x_3 & 0 \\
2x_1 & x_2 & 0 & -x_4 \\
x_3 & x_4 & 0 & 0
\end{vmatrix}
\cdot x_5\cdots x_n.
\]

\subsection{Classification up to dimension $4$}\setcounter{secnumdepth}{5}

We consider the situation of Theorem~\ref{51} and abbreviate $x,y,z,w=x_1,x_2,x_3,x_4$.
By the results of Section~\ref{53}, we may assume that $s=1$ and $n=4$.
Let us first assume that $\Ga(\CC^n,\Der(-\log D))_0$ is a nonsolvable Lie algebra and hence $\langle\sigma_1,\nu_1,\nu_2\rangle=\sl_2$.

Recall that by \cite[IV.4]{Ser87}, $\CC^4$ is a direct sum of irreducible $\sl_2$-modules $W_m$ of dimension $m+1$ and that $W_m$ is represented in a basis $e_0,\dots,e_m$ by
\[
\sigma_1(e_i)=(-m+2i)e_i,\quad
\nu_1(e_i)=(i+1)e_{i+1},\quad
\nu_2(e_i)=(m-i+1)e_{i-1}.
\]
So there are 3 types of $\sl_2$-representations.
The first two cases are $\CC^4=W_2\oplus W_0$ and $\CC^4=W_1\oplus W_1$, which lead to a zero and a nonreduced determinant of the form $\Delta\equiv(xw-yz)^2$ respectively.
But $W_3$ gives the nontrivial linear free divisor
\[
\Delta=\begin{vmatrix}
x&y&z&w\\
-3x&-y&z&3w\\
y&2z&3w&0\\
0&3x&2y&z
\end{vmatrix}
\equiv y^2z^2-4xz^3-4y^3w+18xyzw-27x^2w^2
\]
isomorphic to the discriminant in the space of binary cubics described in Example~\ref{1}.\eqref{1b}.

Now, assume that $\Ga(\CC^n,\Der(-\log D))_0$ is a solvable Lie algebra.
Then, by Lie's Theorem \cite[I.7]{Ser87}, $\nu_1$ and $\nu_2$ can be chosen triangular and also $[\nu_1,\nu_2]$ is triangular.
Hence, $[\nu_1,\nu_2]\in\langle\nu_1,\nu_2\rangle$ and even $[\nu_1,\nu_2]=0$ by nilpotency of $\ad_{\nu_1}$ and $\ad_{\nu_2}$.

\begin{annalif}

We set
\[
\sigma_1=ax\p_x+by\p_y+cz\p_z+dw\p_w
\]
and start a case by case discussion with respect to the cardinality of $\{a,b,c,d\}$. 
In the following we shall omit the details that can be found in the preprint version \cite{GMNS06} of this article.

\subsubsection{$2\leq |\{a,b,c,d\}|\leq 3$}

In each case we refine to subcases depending on whether $\sigma_1(\nu_i)/\nu_i$, $i=1,2$, is zero or not.
All these subcases lead to $\Delta=0$ or $\Delta$ being divisible by a square of a variable, in contradiction to our assumptions. 

\subsubsection{$|\{a,b,c,d\}|=4$}

Since $\nu_1$ and $\nu_2$ are $\sigma_1$-homogeneous and might be chosen triangular, $\sigma(\nu_1)/\nu_1\ne0\ne\sigma(\nu_2)/\nu_2$ and hence $\nu_1$ and $\nu_2$ have at most $3$ nonzero entries by Lemma~\ref{52}.
Using, if necessary, permutations of the basis vectors, we have only to consider the following two cases: 

\paragraph{$\nu_1$ has one nonzero term or $\nu_1$ and $\nu_2$ have at most two nonzero terms.}

In both cases, a detailed discussion leads to the contradiction that $\Delta$ is divisible by a square of a variable.

\paragraph{$\nu _1 $ has three nonzero terms.}

This turns out to be the only case that leads to a linear free divisor with solvable Lie algebra and $s=1$.
We may assume that
\[
\nu_1=
\begin{pmatrix}
0 & 1 & 0 & 0 \\
0 & 0 & 1 & 0 \\
0 & 0 & 0 & 1 \\
0 & 0 & 0 & 0 
\end{pmatrix}
\]
which implies also that $(a,b,c,d)=a\cdot(1,1,1,1)+(0,\lambda,2\lambda,
3\lambda)$ for some $0\neq\lambda\in\QQ$.
The relation $[\nu_1,\nu_2]=0$ then implies that 
\[
\nu_2=p\cdot\nu_1+
\begin{pmatrix}
0 & 0 & q & r \\
0 & 0 & 0 & q \\
0 & 0 & 0 & 0 \\
0 & 0 & 0 & 0 
\end{pmatrix}.
\] 
So using the $\sigma_1$-homogeneity of $\nu_2$ the only case which was not yet considered is
\[
\nu_2= 
\begin{pmatrix}
0 & 0 & 1 & 0 \\
0 & 0 & 0 & 1 \\
0 & 0 & 0 & 0 \\
0 & 0 & 0 & 0
\end{pmatrix},\quad
\Delta\equiv
\begin{vmatrix}
x & y & z & w \\
0 & \lambda y & 2\lambda z & 3\lambda w \\
0 & x & y & z \\
0 & 0 & x & y
\end{vmatrix}
\equiv x(y^3-3xyz+3x^2w).
\] 
The trace equation in Theorem~\ref{51}.\ref{51d} for $\sigma_1$ reads $-\lambda-2\lambda+4a+6\lambda=0$ or $a=-\frac{3}{4}\lambda $.
Setting $\lambda=4$, we obtain $\sigma_1=-3x\p_x+y\p_y+5z\p_z+9w\p_w$.

\end{annalif}
\begin{arxiv}

\subsubsection{$\sigma_1=ax\p_x+ay\p_y+az\p_z+bw\p_w$, $a\ne b$}

Then $f=0$ if $\nu_1$ or $\nu_2$ has only a $\p_w$-component.
We shall tacitly omit this case in the following.

\paragraph{$\sigma_1(\nu_1)/\nu_1=0=\sigma_1(\nu_2)/\nu_2$.}

Then 
\[
\nu_1=
\begin{pmatrix}
0 & 1 & 0 & 0 \\
0 & 0 & \epsilon & 0 \\
0 & 0 & 0 & 0 \\
0 & 0 & 0 & 0 
\end{pmatrix},\quad
\nu_2=
\begin{pmatrix}
0 & * & * & 0 \\
0 & 0 & * & 0 \\
0 & \overline\epsilon & 0 & 0 \\
0 & 0 & 0 & 0 
\end{pmatrix},\quad
f\equiv
\begin{vmatrix}
x & y & z & 0 \\
0 & 0 & 0 & w \\
0 & x & \epsilon y & 0 \\
0 & \overline\epsilon z & * & 0 
\end{vmatrix}
\]
where $\epsilon\in\{0,1\}$ and $\overline\epsilon=1-\epsilon$.
Hence, $f$ is divisible by $x^2$ and not reduced.

\paragraph{$\sigma_1(\nu_1)/\nu_1=0\ne\sigma_1(\nu_2)/\nu_2$.}

Then, after a linear coordinate change in $x,y,z$, we have
\[
\nu_2=
\begin{pmatrix}
0 & 0 & 0 & 0 \\
0 & 0 & 0 & 0 \\
0 & 0 & 0 & 0 \\
1 & 0 & 0 & 0 
\end{pmatrix},\quad
f\equiv
\begin{vmatrix}
x & y & z & 0 \\
0 & 0 & 0 & w \\
* & * & * & * \\
w & 0 & 0 & 0
\end{vmatrix}.
\]
Hence, $f$ is divisible by $w^2$ and not reduced.

\paragraph{$\sigma_1(\nu_1)/\nu_1\ne0\ne\sigma_1(\nu_2)/\nu_2$.}

Then 
\[
\nu_1=
\begin{pmatrix}
0 & 0 & 0 & 0 \\
0 & 0 & 0 & 0 \\
0 & 0 & 0 & 0 \\
1 & 0 & 0 & 0 
\end{pmatrix},\quad
\nu_2=
\begin{pmatrix}
0 & 0 & 0 & 0 \\
0 & 0 & 0 & 0 \\
0 & 0 & 0 & 0 \\
r & s & t & 0 
\end{pmatrix}
\]
and hence $f$ is divisible by $w^3$ and not reduced.

\subsubsection{$\sigma_1=ax\p_x+ay\p_y+bz\p_z+bw\p_w$, $a\ne b$}

Then $f=0$ if $\nu_1$ and $\nu_2$ have only a $\p_z$- and $\p_w$-component.
We shall tacitly omit this case in the following.

\paragraph{$\sigma_1(\nu_1)/\nu_1=0=\sigma_1(\nu_2)/\nu_2$.}

Then 
\[
\nu_1=
\begin{pmatrix}
0 & 1 & 0 & 0 \\
0 & 0 & 0 & 0 \\
0 & 0 & 0 & 0 \\
0 & 0 & 0 & 0 
\end{pmatrix},\quad
\nu_2=
\begin{pmatrix}
0 & 0 & 0 & 0 \\
0 & 0 & 0 & 0 \\
0 & 0 & 0 & 1 \\
0 & 0 & 0 & 0 
\end{pmatrix},\quad
f\equiv
\begin{vmatrix}
x & y & 0 & 0 \\
0 & 0 & z & w \\
0 & x & 0 & 0 \\
0 & 0 & 0 & z
\end{vmatrix}\equiv x^2z^2.
\]

\paragraph{$\sigma_1(\nu_1)/\nu_1=0\ne\sigma_1(\nu_2)/\nu_2$.}

Then 
\[
\nu_1=
\begin{pmatrix}
0 & 1 & 0 & 0 \\
0 & 0 & 0 & 0 \\
0 & 0 & 0 & \epsilon \\
0 & 0 & 0 & 0 
\end{pmatrix},\quad
\nu_2=
\begin{pmatrix}
0 & 0 & * & * \\
0 & 0 & * & * \\
0 & 0 & 0 & 0 \\
0 & 0 & 0 & 0 
\end{pmatrix},\quad
f\equiv
\begin{vmatrix}
x & y & 0 & 0 \\
0 & 0 & z & w \\
0 & x & 0 & \epsilon z \\
0 & 0 & * & *
\end{vmatrix}
\]
where $\epsilon\in\{0,1\}$. 
Hence, $f$ is divisible by $x^2$ (or $z^2$ for transposed $\nu_2$) and 
not reduced.

\paragraph{$\sigma_1(\nu_1)/\nu_1\ne0\ne\sigma_1(\nu_2)/\nu_2$.}

Then 
\[
\nu_1=
\begin{pmatrix}
0 & 0 & 1 & 0 \\
0 & 0 & 0 & 0 \\
0 & 0 & 0 & 0 \\
0 & 0 & 0 & 0 
\end{pmatrix},\quad
\nu_2=
\begin{pmatrix}
0 & 0 & 0 & 0 \\
0 & 0 & 0 & 0 \\
0 & 0 & 0 & 0 \\
0 & 1 & 0 & 0 
\end{pmatrix},\quad
f\equiv
\begin{vmatrix}
x & y & 0 & 0 \\
0 & 0 & z & w \\
0 & 0 & x & 0 \\
0 & w & 0 & 0
\end{vmatrix}\equiv x^2w^2.
\]

\subsubsection{$\sigma_1=ax\p_x+ay\p_y+bz\p_z+cw\p_w$, $|\{a,b,c\}|=3$}

Then $f=0$ if $\nu_1$ and $\nu_2$ have only a $\p_z$- and $\p_w$-component.
We shall tacitly omit this case in the following.

\paragraph{$\sigma_1(\nu_1)/\nu_1=0=\sigma_1(\nu_2)/\nu_2$.}

Then $\nu_1$ and $\nu_2$ are linearly dependent.

\paragraph{$\sigma_1(\nu_1)/\nu_1=0\ne\sigma_1(\nu_2)/\nu_2$.}

Then 
\[
\nu_1=
\begin{pmatrix}
0 & 1 & 0 & 0 \\
0 & 0 & 0 & 0 \\
0 & 0 & 0 & 0 \\
0 & 0 & 0 & 0 
\end{pmatrix},\quad
\nu_2=
\begin{pmatrix}
0 & 0 & * & * \\
0 & 0 & 0 & 0 \\
0 & * & 0 & * \\
0 & * & * & 0 
\end{pmatrix},\quad
f\equiv
\begin{vmatrix}
x & y & z & w \\
0 & 0 & bz & cw \\
0 & x & 0 & 0 \\
0 & * & * & *
\end{vmatrix}.
\]
Hence, $f$ is divisible by $x^2$ and not reduced.

\paragraph{$\sigma_1(\nu_1)/\nu_1\ne0\ne\sigma_1(\nu_2)/\nu_2$.}

Then there are 3 double cases:
\begin{enumerate}

\item
\[
\nu_1=
\begin{pmatrix}
0 & 0 & 1 & 0 \\
0 & 0 & 0 & 0 \\
0 & 0 & 0 & 0 \\
0 & 0 & 0 & 0 
\end{pmatrix},\quad
\nu_2=
\begin{pmatrix}
0 & 0 & * & * \\
0 & 0 & * & * \\
0 & 0 & 0 & 0 \\
0 & * & * & 0 
\end{pmatrix},\quad
f\equiv
\begin{vmatrix}
x & y & z & w \\
0 & 0 & bz & cw \\
0 & 0 & x & 0 \\
0 & * & * & *
\end{vmatrix}.
\]
Hence, $f$ is divisible by $x^2$ (or $z^2$ for transposed $\nu_1$ and 
$\nu_2$) and not reduced.

\item
\[
\nu_1=
\begin{pmatrix}
0 & 0 & 0 & 0 \\
0 & 0 & 0 & 0 \\
0 & 0 & 0 & 1 \\
0 & 0 & 0 & 0 
\end{pmatrix},\quad
\nu_2=
\begin{pmatrix}
0 & 0 & 0 & * \\
0 & 0 & 0 & * \\
* & * & 0 & 0 \\
0 & 0 & 0 & 0 
\end{pmatrix},\quad
f\equiv
\begin{vmatrix}
x & y & z & w \\
0 & 0 & bz & cw \\
0 & 0 & 0 & z \\
* & * & 0 & *
\end{vmatrix}.
\]
Hence, $f$ is divisible by $z^2$ (or $w^2$ for transposed $\nu_1$ and 
$\nu_2$) and not reduced.

\item
\[
\nu_1=
\begin{pmatrix}
0 & 0 & 1 & 0 \\
0 & 0 & 0 & 0 \\
0 & 0 & 0 & 1 \\
0 & 0 & 0 & 0 
\end{pmatrix},\quad
\nu_2=
\begin{pmatrix}
0 & 0 & 0 & * \\
0 & 0 & 0 & * \\
0 & 0 & 0 & 0 \\
0 & 0 & 0 & 0 
\end{pmatrix},\quad
f\equiv
\begin{vmatrix}
x & y & z & w \\
0 & 0 & bz & cw \\
0 & 0 & x & z \\
0 & 0 & 0 & *
\end{vmatrix}=0
\]
or $f$ is divisible by $w^2$ for transposed $\nu_1$ and $\nu_2$ and not 
reduced.

\end{enumerate}

\subsubsection{$\sigma_1=ax\p_x+by\p_y+cz\p_z+dw\p_w$, $|\{a,b,c,d\}|=4$}

Since $\nu_1$ and $\nu_2$ are $\sigma_1$-homogeneous and might be chosen 
triangular, $\sigma(\nu_1)/\nu_1\ne0\ne\sigma(\nu_2)/\nu_2$ and hence 
$\nu_1$ and $\nu_2$ have at most $3$ non-zero entries by Lemma~\ref{52}.
Using if necessary permutations of the basis vectors, we have only to 
consider the following cases: 

\paragraph{$\nu _1$ has one non-zero term.}

We may assume that 
\[
\nu_1=
\begin{pmatrix}
0 & 1 & 0 & 0 \\
0 & 0 & 0 & 0 \\
0 & 0 & 0 & 0 \\
0 & 0 & 0 & 0 
\end{pmatrix}.
\] 
Since $[\nu_1,\nu_2]=0$, we find that the first column of $\nu_2$ is 
zero and hence
\[
f\equiv
\begin{vmatrix}
x & y & z & w \\
ax & by & cz & dw \\
0 & x & 0 & 0 \\
0 & * & * & *
\end{vmatrix}
\]
is non-reduced as a multiple of $x^2$. 
In what follows, we may also exclude all cases where $\nu_2$ has only 
one term by exchanging the roles of $\nu_1$ and $\nu_2$. 

\paragraph{$\nu_1$ has two non-zero terms.}
We have two cases:
\[
\nu_1=
\begin{pmatrix}
0 & 1 & 0 & 0 \\
0 & 0 & 1 & 0 \\
0 & 0 & 0 & 0 \\
0 & 0 & 0 & 0 
\end{pmatrix}\quad
\text{ or }\quad
\nu_1=
\begin{pmatrix}
0 & 1 & 0 & 0 \\
0 & 0 & 0 & 0 \\
0 & 0 & 0 & 1 \\
0 & 0 & 0 & 0 
\end{pmatrix}.
\] 
In the first case, since $[\nu_1,\nu_2]=0$, we obtain that $\nu _2$ has 
the form
\[
\nu_2=p\cdot\nu_1+
\begin{pmatrix}
0 & 0 & q & r \\
0 & 0 & 0 & 0 \\
0 & 0 & 0 & 0 \\
0 & 0 & s & 0 
\end{pmatrix}.
\] 
Since $a-b=b-c$ because of the homogeneity of $\nu_1$, the $\sigma_1$-degrees 
corresponding to $r,s$ are $a-d\neq d-c$.
Then, because of the homogeneity of $\nu_2$, we are in the situation 
where $\nu _2$ has only one term. 
In the second case, we obtain 
\[
\nu_2=p\cdot\nu_1+
\begin{pmatrix}
0 & 0 & u & v \\
0 & 0 & 0 & u \\
r & s & 0 & q \\
0 & r & 0 & 0 
\end{pmatrix}.
\] 
If $u=1$ and $v=r=s=q=0$ then 
\[
\nu_2=p\cdot\nu_1+
\begin{pmatrix}
0 & 0 & 1 & 0 \\
0 & 0 & 0 & 1 \\
0 & 0 & 0 & 0 \\
0 & 0 & 0 & 0 
\end{pmatrix},\quad
f\equiv
\begin{vmatrix}
x & y & z & w \\
ax & by & cz & dw \\
0 & x & 0 & z \\
0 & 0 & x & y
\end{vmatrix}.
\]
So $f$ is non-reduced as a multiple of $x^2$.  
Similarly, if $r=1$ and $u=v=s=q=0$ then $f$ is a multiple of $z^2$.
In the other cases,  either $\nu _2$ has only one non-zero term or 
three non-zero terms ($u=s=1$).
This latter case we shall study now.

\paragraph{$\nu _1 $ has three non-zero terms.}
We may assume that
\[
\nu_1=
\begin{pmatrix}
0 & 1 & 0 & 0 \\
0 & 0 & 1 & 0 \\
0 & 0 & 0 & 1 \\
0 & 0 & 0 & 0 
\end{pmatrix}
\]
which implies also that $(a,b,c,d)=a\cdot(1,1,1,1)+(0,\lambda,2\lambda,
3\lambda)$ for some $0\neq\lambda\in\QQ$.
The relation $[\nu_1,\nu_2]=0$ then implies that 
\[
\nu_2=p\cdot\nu_1+
\begin{pmatrix}
0 & 0 & q & r \\
0 & 0 & 0 & q \\
0 & 0 & 0 & 0 \\
0 & 0 & 0 & 0 
\end{pmatrix}.
\] 
So the only remaining possibility not yet considered is
\[
\nu_2-p\nu _1= 
\begin{pmatrix}
0 & 0 & 1 & 0 \\
0 & 0 & 0 & 1 \\
0 & 0 & 0 & 0 \\
0 & 0 & 0 & 0
\end{pmatrix},\quad
f\equiv
\begin{vmatrix}
x & y & z & w \\
0 & \lambda y & 2\lambda z & 3\lambda w \\
0 & x & y & z \\
0 & 0 & x & y
\end{vmatrix}
\equiv x(y^3-3xyz+3x^2w).
\] 
The equation $\sigma_1(\nu_1)/\nu_1+\sigma_1(\nu_2)/\nu_2+\tr(\sigma_1)=0$ 
becomes $-\lambda-2\lambda+4a+6\lambda=0$ or $a=-\frac{3}{4}\lambda $.
Setting $\lambda=4$, we obtain $\sigma_1=-3x\p_x+y\p_y+5z
\p_z+9w\p_w$ and $\sigma_1(f)=0$.

\end{arxiv}

\subsection{Summary of the classification up to dimension 4}\label{56}

The following Table~\ref{88} summarizes our classification of linear free divisors up to dimension $4$.
The matrices are interpreted row-wise as bases of $\Der(-\log D)$.

\begin{longtable}{|c|c|c|c|c|}
\caption{Classification of linear free divisors up to dimension $4$}\label{88}\\
\hline
\rule{0pt}{11pt}$n$ & $\Delta$ & $\Der(-\log D)$ & $\gg_D$ & reductive?\\
\hline
\hline
$1$ & $x$&$\begin{pmatrix}x\end{pmatrix}$ & $\CC$ & Yes \\
\hline
\hline
$2$ & $xy$&$\begin{pmatrix}x&0\\0&y\end{pmatrix}$ & $\CC^2$ & Yes \\
\hline
\hline
$3$ & $xyz$&$\begin{pmatrix}x&0&0\\0&y&0\\0&0&z\end{pmatrix}$ & $\CC^3$ & Yes \\
\hline
$3$& $(y^2+xz)z$&$\begin{pmatrix}x&y&z\\4x&y&-2z\\-2y&z&0\end{pmatrix}$ & $\bb_2$ & No \\
\hline
\hline
$4$ & $xyzw$&$\begin{pmatrix}
x&0&0&0\\0&y&0&0\\0&0&z&0\\0&0&0&w\end{pmatrix}$ & $\CC^4$ & Yes \\
\hline
$4$ & $(y^2+xz)zw$&$\begin{pmatrix}
x&y&z&0\\4x&y&-2z&0\\-2y&z&0&0\\0&0&0&w\end{pmatrix}$ & $\CC\oplus\bb_2$ & No \\
\hline
$4$ & $(yz+xw)zw$&
$\begin{pmatrix}
x&0&0&-w\\0&y&0&w\\0&0&z&w\\z&-w&0&0
\end{pmatrix}$ 
& $\CC^2\oplus \gg_0$ & No \\
\hline
$4$ & $x(y^3-3xyz+3x^2w)$&
$\begin{pmatrix}
x&y&z&w\\0&y&2z&3w\\
0&x&y&z\\0&0&x&y
\end{pmatrix}$ 
& $\CC\oplus \gg$ & No \\
\hline
$4$ & \begin{minipage}{3cm}\begin{center}$y^2z^2-4xz^3-4y^3w+18xyzw-27w^2x^2$\end{center}\end{minipage}&
$\begin{pmatrix}
3x&2y&z&0\\0&3x&2y&z\\y&2z&3w&0\\0&y&2z&3w
\end{pmatrix}$ & $\gl_2(\CC)$ & Yes \\
\hline
\end{longtable}

The annihilators of $\Delta$ in the Lie algebras for $\Delta=(yz+xw)zw$ and $\Delta=x(y^3-3xyz+3x^2w)$ are described in \cite[Ch.~I, \S4]{Jac79}.
The former is the direct sum of $\CC$ and the non-Abelian Lie algebra $\gg_0$ of dimension $2$, and the latter is the $3$-dimensional Lie algebra $\gg$ characterized as having $2$-dimensional Abelian derived algebra $\gg'$, on which the adjoint action of a basis vector outside $\gg'$ is semi-simple with eigenvalues $1$ and $2$. 
Straightforward computations show that the two \emph{groups} $G_D^\circ$ are, respectively, the set of $4\times 4$ matrices of the form
\[
\begin{pmatrix}
x^{-1}y^{-2}&0&z&0\\
0&x^{-2}y^{-1}&0&-x^{-1}yz\\
0&0&x&0\\
0&0&0&y
\end{pmatrix}
\quad
\mbox{and}
\quad
\begin{pmatrix}
x^{-3}&0&0&0\\
y&x&0&0\\
z&x^4y&x^5&0\\
x^3yz-\frac{1}{3}x^6y^3&x^4z&x^8y&x^9
\end{pmatrix}
\] 
with $x,y\in\CC^*$ and $z\in \CC$ in the first and $x\in\CC^*, y,z\in\CC$ in the second.

\section{Strong Euler homogeneity and local quasihomogeneity}\label{5}

In this section, we investigate linear free divisors with respect to the properties of local quasihomogeneity and strong Euler homogeneity from Definitions~\ref{80} and \ref{81}.

The following reformulation of the definition of local quasihomogeneity is a direct consequence of the Poincar\'e--Dulac Theorem \cite[Ch.~3, \S 3.2]{ABGI88} and  Artin's Approximation Theorem \cite{Art68}.

\begin{theo}
A divisor $D$ is locally quasihomogeneous if and only if, at any $p\in D$, there is an Euler vector field $\chi$ for $D$ at $p$ whose degree zero part $\chi_0$ has strictly positive eigenvalues.
\end{theo}

We denote by $\D=\D_{\CC^n}$ the sheaf of germs of linear differential operators with holomorphic coefficients on $\CC^n$.
It is naturally equipped with an increasing filtration $F$ of coherent $\O$-modules by the order of differential operators and we denote by $\sigma(P)$ the symbol of $P\in\D$ in $\gr_F\D$.
Note that $\gr_F\D_p\cong\O_p[\p]=\CC\{x\}[\p]$ in a local coordinate system $x$ at $p$, where we identify $\sigma(\p_i)=\p_i$.The following property is closely related to local quasihomogeneity.

\begin{defi}
A free divisor $D$ is called \emph{Koszul free} if, at any $p\in D$, there exists a basis $\delta_1,\dots,\delta_n$ of $\Der(-\log D)_p$ such that $\sigma(\delta_1),\dots,\sigma(\delta_n)$ is a regular sequence in $\gr_F\D_p$.
\end{defi}

Koszul freeness can be interpreted geometrically in terms of the \emph{logarithmic stratification}, introduced by K.~Saito \cite{Sai80}, which is the partition of $D$ into the integral varieties of the distribution $\Der(-\log D)$. 
As it is not always locally finite, the term ``stratification'' is a misnomer, but is generally used. 
If $D_\alpha$ is a stratum of the logarithmic stratification and $p\in D_\alpha$ then $T_pD_\alpha=\Der(-\log D)(p)$.
The graded ring $\gr_F\D_p\cong\O_p[\p]$ contains $\Der_p=\oplus_{i=1}^n\O_p\p_i$ and can be identified with the ring of functions on the cotangent space $T^*(\CC^n,p)$ of the germ $(\CC^n,p)$, polynomial on the fibers and analytic on the base.

\begin{defi}
The \emph{logarithmic characteristic variety} $L_{\CC^n}(D)$ of $D$ is the variety in $T^*\CC^n$ defined by the image of $\Der(-\log D)$ in $\gr_F\D$.
\end{defi}

Thus $D$ is Koszul free at $p$ if and only if $L_{\CC^n}(D)$ is purely $n$-dimensional \cite[1.8]{CN02}. 
Moreover, by \cite[3.16]{Sai80}, $L_{\CC^n}(D)$ is the union, over all strata $D_\alpha\subseteq D$ in the logarithmic stratification of $D$, of the conormal bundle $T^*_{D_\alpha}\CC^n$ of $D_\alpha$, each of which is $n$-dimensional.
This proves the following result.

\begin{theo}
A free divisor $D$ is Koszul free if and only if the logarithmic stratification is locally finite. 
\end{theo}

A locally quasihomogeneous free divisor is Koszul free \cite[4.3]{CN02}, and thus the failure of Koszul-freeness serves as a computable criterion for the failure of local quasihomogeneity.

As in Section \ref{49}, let $D\subseteq\CC^n$ be a linear free divisor and $\chi,\delta_2,\dots,\delta_n$ a global degree $0$ basis of $\Der(-\log D)$ with $\delta_i(\Delta)=0$ for $\Delta=\det(S)$ where
\[
S:=\begin{pmatrix}\chi(x_j)\\\delta_i(x_j)\end{pmatrix}_{i,j}
\]
The arguments which follow are valid as well for an arbitrary germ of a free divisor $D\subseteq(\CC^n,0)$ with a germ of an Euler vector field $\chi\in\Der_0$ at $0$. 

The following criterion gives a method to test strong Euler homogeneity algorithmically.
The reduced variety $S_k$ defined by the $(k+1)\times(k+1)$-minors of the $n\times n$-matrix $S$ is the union of logarithmic strata of dimension at most $k$.
In more invariant terms, $S_k$ is the variety of zeros of the $(n-k-1)$'st Fitting ideal of the $\O_{\CC^n}$-module $\Der/\Der(-\log D)$.
Thus a free divisor $D$ is Koszul free if and only if $\dim S_k\leq k$ for all $k$. 
Note that $S_n=\CC^n$, $S_{n-1}=D$, and $S_{n-2}=\Sing(D)$.
For a linear free divisor, $S_0=\{0\}$ because of the presence of the Euler vector field. 
Since $\dim\Sing(D)<\dim D$, it follows that linear free divisors are Koszul free in 
dimension $n\le3$.
 
In order to characterize strong Euler homogeneity, we also consider the reduced variety $T_k\supseteq S_k$ defined by the $(k+1)\times(k+1)$-minors of the $(n-1)\times n$-matrix 
\[
T:=(\delta_i(x_j))_{i,j}.
\]
Again, this is the variety defined by a Fitting ideal, this time the $(n-k-2)$nd Fitting ideal of the module $\Der/\Der(-\log\Delta)$.
Note that, by definition,
\[
S=\begin{pmatrix}x\\T\end{pmatrix}.
\]

\begin{lemm}\label{57}
$D$ is strongly Euler homogeneous if and only if $S_k=T_k$ for $0\le k\le n-2$.
\end{lemm}

\begin{proof}
A vector field $\delta\in\big(\mm_p\cdot\chi+\sum_i\O_p\cdot\delta_i\big)\cap\mm_p\cdot\Der_p$ is not an Euler vector field at $p\in D$ since $\delta(\Delta)\in\mm_p\cdot\Delta$, and indeed $\delta(u\cdot\Delta)\in\mm_p\cdot u\Delta$ for any unit $u$.
Hence, an Euler vector field $\eta$ for $D$ at $p$ must be of the form $\eta=a_0\cdot\chi+\sum_ia_i\cdot\delta_i\in\mm_p\cdot\Der_p$ with $a_0(0)\neq 0$. 
This means that $\chi(p)\in\sum_i\CC\cdot\delta_i(p)$, and the matrices $S$ and $T$ have equal rank at $p$. 
Conversely if $\chi(p)=\sum \lambda_i\delta_i(p)$ then up to multiplication by a scalar, $\chi-\sum_i\lambda_i\delta_i$ is an Euler field for $D$ at $p$.
\end{proof}

\begin{rema}
The proof of Lemma~\ref{57} shows that the question of local 
quasihomogeneity is much more complicated:
The degree zero parts of Euler vector fields at a point $p\in D$ are the degree zero parts of vector fields $a_1\cdot\chi+\sum_ia_i\cdot\delta_i$ where $a_1,\dots,a_n$ are linear forms such that $a_1(p)\cdot\chi(p)+\sum_ia_i(p)\cdot\delta_i(p)=0$.
\end{rema}

For $k=1,\dots,n$, let $M_k=(-1)^{k+1}\det (\delta_i(x_j))_{j\neq k}$.

\begin{lemm}\label{58}
For $k=1,\dots,n$, $\p_k(\Delta)=n\cdot M_k$.
In particular, $S_k=T_k$ for $k=n-2$.
\end{lemm}

\begin{proof}
Since
\[
S\begin{pmatrix}\p_1(\Delta)\\\vdots\\\p_n(\Delta)\end{pmatrix}=
\begin{pmatrix}\chi\\\delta_2\\\vdots\\\delta_n\end{pmatrix}(\Delta)=
\begin{pmatrix}n\cdot\Delta\\0\\\vdots\\0\end{pmatrix}
\]
we obtain, by canceling $\Delta$,
\[
\begin{pmatrix}\p_1(\Delta)\\\vdots\\\p_n(\Delta)\end{pmatrix}=\check 
S\begin{pmatrix}n\\0\\\vdots\\0\end{pmatrix}=n\cdot\begin{pmatrix}M_1\\
\vdots\\M_n\end{pmatrix}
\]
where $\check S$ denotes the cofactor matrix of $S$.
\end{proof}

\begin{lemm}\label{59}
$S_0=T_0$.
\end{lemm}

\begin{proof}
Assume that $T_0\ne S_0=\{0\}$.
By homogeneity, $T_0$ contains the $x_n$-axis after an appropriate linear 
coordinate change.
Then $T$ is independent of $x_n$.
Writing $x'=x_1,\dots,x_{n-1}$, we have that $\Delta=g+x_n\cdot\Delta'$ where $g$ and $\Delta':=M_n$ depend only on $x'$.
Since $\Delta$ does not depend on fewer variables, we must have $\Delta'\ne0$.
For $i=2,\dots,n$, let $\delta_i':=\sum_{j=1}^{n-1}\delta_i(x_j)\p_j$ 
be the projection of $\delta_i$ to the $\CC[x]$-module with basis 
$\p':=\p_1,\dots,\p_{n-1}$.
Then, for $i=2,\dots,n$, $\delta'_i(\Delta')=0$, as it is the coefficient of $x_n$ in $\delta_i(\Delta)=0$.
Since the rank of the $\CC[x']$-annihilator of $\p_1(\Delta'),\dots,\p_{n-1}(\Delta')$ is strictly smaller than $n-1$, 
there must be a relation 
$\sum_{i=2}^na_i\delta'_i=0$ for some homogeneous polynomials 
$a_i\in\CC[x']$.
But since $\delta_2,\dots,\delta_n$ are independent over $\CC[x]$, $\sum_{i=2}^na_i\delta_i(x_n)\ne0$ and hence $0=\sum_{i=2}^na_i\delta_i(\Delta)=\bigl(\sum_{i=2}^na_i\delta_i(x_n)\bigr)\cdot\Delta'$, contradicting the fact that $\Delta'\ne0$.
\end{proof}

\begin{lemm}\label{60}
Let $D$ be strongly Euler homogeneous.
Then $D$ is locally quasihomogeneous on the complement of $S_{n-3}$.
In particular, $D$ is locally quasihomogeneous if $S_{n-3}=\{0\}$.
\end{lemm}

\begin{proof}
By \cite[3.5]{Sai80}, $(D,p)=(D',p')\times(\CC^{n-2},0)$ for $p\in S_{n-2}\ssm S_{n-3}$ where $(D',p')\subseteq(\CC^2,0)$ is strongly Euler homogeneous by \cite[3.2]{GS06}.
As the germ of a curve, $(D',p')$ has an isolated singularity.
Then $(D',p')$, and hence $(D,p)$, are quasihomogeneous, by Saito's theorem \cite{Sai71}.
\end{proof}

\begin{theo}\label{61}
Every linear free divisor in dimension $n\le 4$ is locally quasihomogeneous and hence LCT and GLCT hold.
\end{theo}

\begin{proof}
By Lemmas \ref{58} and \ref{59}, $S_1=T_1$ if $n=3$ and $S_0=T_0$. 
If $n\le3$ then $D$ is strongly Euler homogeneous by Lemma~\ref{57} and so locally quasihomogeneous by Lemma~\ref{60}.

For $n=4$ analogous arguments yield $S_0=T_0=\{0\}$ and $S_2=T_2$.
Now we use the classification in Subsection \ref{56} and a case by case study:
In each case, one can verify that $S_1=T_1$ and construct an Euler vector field with positive eigenvalues in degree zero at each point of $S_1\smallsetminus S_0$.
Again this is sufficient for local quasihomogeneity by Lemma \ref{60}.
For $\Delta=(yz+xw)zw$, $S_1=\{xy=z=w=0\}$ and $2\chi-\sigma+\frac{x-\xi}{\xi}\sigma$, where $\sigma=2x\p_x+y\p_y-w\p_w$, is an Euler vector field at $(\xi,0,0,0)\not\in S_0$ with eigenvalues $2,1,2,3$ in degree zero.
For $\Delta=x(y^3-3xyz+3x^2w)$, $S_1=\{x=y=z=0\}$ and $9\chi-\sigma+\frac{w-\om}{\om}\sigma$, where $\sigma=-3x\p_x+y\p_y+5z\p_z+9w\p_w$, is an Euler vector field at $(0,0,0,\om)\not\in S_0$ with eigenvalues $12,8,4,9$ in degree zero.
The remaining cases are trivial.

By \cite{CMN96}, local quasihomogeneity implies that LCT and hence, by taking global sections, GLCT holds.
\end{proof}

\subsection{Example~\ref{41} again}\label{62}

In this subsection, we study the linear free divisor in Example~\ref{41} in detail and show that if $n>2$ it is not Koszul free and hence not locally quasihomogeneous.
However we will see that it is strongly Euler homogeneous, like all other linear free divisors whose strong Euler homogeneity has been investigated.

Denote by $x_{i,j}$, $1\leq  i\leq n$, $1\leq  j\leq n+1$, the coordinates on the space of $n\times(n+1)$-matrices $M_{n,n+1}$. 
The Lie group $G=\Gl_n(\CC)\times \Gl_1(\CC)^{n+1}$ acts on $M_{n,n+1}$ by left matrix multiplication of $\Gl_n(\CC)$ and multiplication of the $j$th factor $\Gl_1(\CC)=\CC^*$ on the $j$th column of members of $M_{n,n+1}$.
By Lemma~\ref{9}, $\Der(-\log D)$ is generated by the infinitesimal action of the Lie algebra of $G$ and hence a basis of $\Der(-\log D)$ is extracted from the set of $n^2+n+1$ vector fields
\begin{align}\label{22}
\xi _{i,j}= \sum _{k=1}^{n+1}x_{i,k}\p_{j,k}, \text{ for } 1\leq i,j \leq n, \\
\nonumber\xi _{i}= \sum _{l=1}^{n}x_{l,i}\p_{l,i}, \text{ for } 1\leq i \leq n+1,
\end{align} 
by omitting one because of the relation
\[
\chi=\sum_{i=1}^n\xi _{i,i} = \sum_{j=1}^{n+1}\xi _{j}
\]
corresponding to the Lie algebra of the kernel of the action.
Note that the vector field $\xi _{i,i}$ resp.\ $\xi _{j}$ is the Euler vector field related to the $i$th row resp.\ to the $j$th column of the general $n\times(n+1)$-matrix and that $\chi$ is the global Euler vector field on $M_{n,n+1}$.

Since the determinant $\Delta _j$ has degree one with respect to each line and to each column except the $j$th column for which the degree is zero, the degree of $\Delta$ equals $n+1$ with respect to a row and $n$ with respect to a column. 
These considerations yield
\[
\xi_{i,i}(\Delta)=(n+1)\Delta,\quad
\xi_{j}(\Delta)=n\Delta,\quad
\xi_{i,j}(\Delta)=0 \text{ for } i\ne j,
\]
and one can easily derive a basis of the vector fields annihilating $\Delta$.

The following lemma is self evident by definition of the action of $G$ on $M_{n,n+1}$ and we shall use it implicitly.
In particular, the rank of a $G$-orbit is well-defined as the rank of any of its elements.

\begin{lemm}\
\begin{asparaenum}
\item Two matrices in $M_{n,n+1}$ having the same row space are in the same $G$-orbit.
Similarly two matrices given by lists of column vectors $A= C_1,\dots,C_{n+1}$ and $A'=C'_1,\dots,C'_{n+1}$ are in the same $G$-orbit if there is a $\lambda _j\in\CC^*$ such that $C'_j=\lambda_jC_j$ for all $j=1,\dots,n+1$.
\item If $A$ and $A'$ are in the same $G$-orbit in $M_{n,n+1}$ then any submatrix of $A$ consisting of columns $C_{i_1}, \dots , C_{i_p}$ has the same rank as the submatrix of $A'$ consisting of the corresponding columns $C'_{i_1}, \dots , C'_{i_p}$ of $A'$. 
In particular $A$ and $A'$ have the same rank.
\end{asparaenum}
\end{lemm}

By the left action of $\Gl_n(\CC)\subseteq G$, any $G$-orbit in rank $r$ 
contains, up to permutation of columns, an element of the form
\beq\label{63}
\begin{pmatrix}
1 & \hdots & 0 & x_{1,r+1} & \hdots & x_{1,n+1} \\
\vdots & \ddots & \vdots & \vdots && \vdots \\
0 & \hdots & 1 & x_{r,r+1} & \hdots & x_{r,n+1} \\
0 & \hdots & 0 & 0 & \hdots & 0 \\
\vdots && \vdots & \vdots && \vdots \\
0 & \hdots & 0 & 0 & \hdots & 0 \\
\end{pmatrix}
\eeq
By using also the action of $G$, we may assume that $x_{i,r+1}\in\{0,1\}$:
If $x_{i,r+1}\ne0$ then one can divide the $i$th row by $x_{i,r+1}$ and multiply the $i$th column by $x_{i,r+1}$.
Thus there is only a finite number of maximal rank $G$-orbits including the generic orbit for which all $x_{i,n+1}$ equal $1$. 

\begin{prop}\label{64}\
\begin{asparaenum}
\item There are only finitely many $G$-orbits in $M_{2,3}$, and the linear free divisor $D\subseteq M_{2,3}$ is locally quasihomogeneous\footnote{See Remark \ref{26} below.}.
\item The number of $G$-orbits in the linear free divisor $D\subseteq M_{3,4}$ is infinite.
In particular, the set of $G$-orbits in $D$ is not locally finite, and $D$ is not Koszul free and hence not locally quasihomogeneous.
\end{asparaenum}
\end{prop}

\begin{proof}\
\begin{asparaenum}
\item The first statement follows, by Gabriel's theorem \cite{Gab72},
from the fact that we are considering the representation 
space of a Dynkin quiver, here of type $D_4$.
In fact, in the case of $M_{2,3}$, the only orbits which 
remain to be considered 
are $\{0\}$ and the rank one orbits which contain, up to permutation of 
columns, one of the typical elements:
\[ 
\begin{pmatrix}
1&1&1\\
0&0&0                       
\end{pmatrix}, 
\begin{pmatrix}
1&1&0\\
0&0&0                       
\end{pmatrix},
\begin{pmatrix}
1&0&0\\
0&0&0                        
\end{pmatrix}.
\] 
At each point $x\neq 0$, $D$ is isomorphic to the product of the germ at $x$ of the orbit $I_x$ of $x$, and the germ at $x$ of $D':=D\cap T$, where $T$ is a smooth transversal to $I_x$ of complementary dimension.
Since $T$ is logarithmically transverse to $D$ in the neighborhood of $x$, $D'$ is a free divisor.
By the Cancellation Property for products of analytic spaces \cite{HM90}, $(D',x)$ is determined up to isomorphism by the fact that $(D,x)\simeq I_x\times(D',x)$, so it does not matter which transversal to $I_x$ we choose. 
In the Table~\ref{90}, we take $T$ to be affine. 
The local equations of $D$ shown in the last column are simply the restriction of the original equation of $D$ to the transversal
$T$. 
By inspection of these equations, $D$ is locally quasihomogeneous.

\begin{longtable}{|c|c|l|}
\caption{Local analysis of $G$-orbits in $M_{2,3}$}\label{90}\\
\hline
\rule{0pt}{11pt}Representative & Transversal & Local Equation\\
\hline
$\begin{pmatrix}1&0&0\\0&0&1\end{pmatrix}$&
$\begin{pmatrix}1&x_{12}&0\\0&x_{22}&1\end{pmatrix}$&
$x_{21}x_{22}=0$\\
\hline
$\begin{pmatrix}1&1&1\\0&0&0\end{pmatrix}$&
$\begin{pmatrix}1&1&1\\x_{21}&x_{22}&0\end{pmatrix}$&
$x_{21}x_{22}(x_{22}-x_{21})=0$\\
\hline
$\begin{pmatrix}1&1&0\\0&0&0\end{pmatrix}$&
$\begin{pmatrix}1&1&x_{13}\\0&x_{22}&x_{23}\end{pmatrix}$&
$x_{22}x_{23}(x_{23}-x_{22}x_{13})=0$\\
\hline
$\begin{pmatrix}1&0&0\\0&0&0\end{pmatrix}$&
$\begin{pmatrix}1&x_{12}&x_{13}\\0&x_{22}&x_{23}\end{pmatrix}$&
$x_{22}x_{23}(x_{12}x_{23}-x_{22}x_{13})=0$\\
\hline
\end{longtable}

\item 
In the case of $M_{3,4}$, consider the stratum in $D$ consisting of matrices of rank 2. 
The four columns span a $2$-dimensional plane, and assuming they are pairwise independent, determine four lines in this plane.
The cross ratio of these four lines is a $G_D$ invariant:
quadruples spanning the same plane, but with different cross-ration, cannot be equivalent. 
Thus there are infinitely many orbits.
Now by \cite[4.3]{CN02} $D$ is not locally quasihomogeneous.
\end{asparaenum}
\end{proof}

\begin{prop}
The linear free divisor $D\subseteq M_{n,n+1}$ from Example~\ref{41} is 
strongly Euler homogeneous for any $n$.
\end{prop}

\begin{proof}
Let us consider a rank $r$ orbit of $G$ in $M_{n,n+1}$.
If $r<n$, we can find a point $A$ in this orbit with a zero row, say row number $i$.
Then the Euler vector field $\xi_{i,i}$ of this row is an Euler vector field at $A$. 

If $r=n$ we can assume that $A$ is of the form \eqref{63} with $x_{i,n+1}=1$ for $1\leq i\leq s$ and 
$x_{i,n+1}=0$ for $s+1\leq i \leq n$ for some $s\leq n$.
Then by \eqref{22} the space parametrized by the variables $x_{i,n+1}$ with $s+1\leq i \leq n$ is a smooth transversal to the orbit at $A$ and the restricted equation of $D$ is just $x_{s+1,n}\cdots x_{n,n+1}=0$.
Thus $D$ is normal crossing and hence strongly Euler homogeneous. 
\end{proof}

\subsection{Example~\ref{40} again}\label{65}

In this subsection, we show that the linear free divisors in Example \ref{40} are locally quasihomogeneous and hence Koszul free by \cite[4.3]{CN02}.
By \cite{CMN96}, this implies that LCT holds although the defining group is not reductive.

We denote by $x_{i,j}$, $1\leq  i\leq j\leq n$, the coordinates on the 
space of symmetric $n\times n$-matrices $\Sym_n(\CC)\subseteq M_{n,n}$. 
Let $D\subseteq\Sym_n(\CC)$ be the divisor defined by the product
\[
\Delta={\det}_1\cdots{\det}_n
\]
of minors 
\[
{\det}_k=
\begin{vmatrix}
x_{1,1} & \cdots & x_{1,k} \\
\vdots & & \vdots \\
x_{k,1} & \cdots & x_{k,k}
\end{vmatrix}.
\]

By Example~\ref{40}, the group $B_n\subseteq\Gl_n(\CC)$ of upper triangular matrices acts on $\Sym_n(\CC)$ by transpose conjugation
\[
B\cdot S=B^tSB,\text{ for }B\in B_n, S\in\Sym_n(\CC)
\]
and the discriminant $D$ is a linear free divisor.
Thus, $\Der(-\log D)$ can be identified with the Lie algebra of $B_n$ and has a basis consisting of the $\frac{1}{2}n(n+1)$ vector fields
\[
\xi_{i,j}=x_{1,i}\del 1 j +\cdots+x_{i,i}\del i j +\cdots+2x_{i,j}\del j j +
\cdots+x_{i,n}\del j n \text{ for }1\leq i\le j\leq n.
\]
It may be helpful to view this as the symmetric matrix
\[
\begin{pmatrix}
0&\cdots&0&\cdots&0&x_{1,i}&0&\cdots&0\\
\vdots&&\vdots&&\vdots&\vdots&\vdots&&\vdots\\
0&\cdots&0&\cdots&0&x_{i,i}&0&\cdots&0\\
\vdots&&\vdots&&\vdots&\vdots&\vdots&&\vdots\\
0&\cdots&0&\cdots&0&x_{i,j-1}&0&\cdots&0\\
x_{1,i}&\cdots&x_{i,i}&\cdots&x_{i,j-1}&2x_{i,j}&x_{i,j+1}&
\cdots&x_{i,n}\\
0&\cdots&0&\cdots&0&x_{i,j+1}&0&\cdots&0\\
\vdots&&\vdots&&\vdots&\vdots&\vdots&&\vdots\\
0&\cdots&0&\cdots&0&x_{i,n}&0&\cdots&0
\end{pmatrix}
\]
in which all the nonzero elements lie in the $j$'th row and the $j$'th column.
Note that the Euler vector field is
\[
\chi=\frac{1}{2}\sum_{i=1}^n\xi _{i,i}.
\]
For $i<j$, $\xi_{i,j}$ is nilpotent, so that $\xi_{i,j}(\Delta)=0$. 
The vector field $\xi_{i,i}$ is the infinitesimal generator of the $\CC^*$ action in which the $i$'th row and column are simultaneously multiplied by $\lambda\in\CC^*$. 
It follows that each determinant $\det_k$ with $k\geq i$ is homogeneous of degree $2$ with respect to $\xi_{i,i}$, and we conclude that
\[
\xi _{i,i}(\Delta)=2(n-i+1)\Delta, \quad \xi _{i,j}(\Delta)=0 \text{ for } i<j.
\]

\begin{lemm}\label{66}
There are finitely many $B_n$-orbits in $\Sym_n(\CC)$. 
\end{lemm}

\begin{proof}
If $i\leq j$, the pair of elementary row and column operations (``add $c$ times column $i$ to column $j$'', ``add $c$ times row $j$ to row $i$'') can be effected by the action of $B_n$. 
By such operations any symmetric matrix may be brought to a normal form with at most a single nonzero element in each row and column. 
Another operation in $B_n$ changes each of these nonzero elements to a $1$.
Thus there are only finitely many $B_n$-orbits in $\Sym_n(\CC)$. 
\end{proof}

By the discussion at the start of Section~\ref{5}, it follows that $D$ is 
Koszul free. In fact this will also follow from  

\begin{prop}\label{67}
The linear free divisor $D$ of Example~\ref{40} associated with the action of  $B_n$ on $\Sym _n(\CC)$ is locally quasi-homogeneous.
\end{prop}

To prove this, it is enough to show that at each point $S$ of $D$ there is an element of $\Der(-\log D)_S$ which vanishes at $S$ and whose linear part is diagonal with positive eigenvalues. 
This is the result of the proposition below.  
In what follows we fix a symmetric matrix $S$ such that $s_{i,j}\in \{ 0,1\}$, with at most one nonzero coefficient in each row and column. 
By Lemma~\ref{66}, each $B_n$ orbit contains such a matrix, and local quasihomogeneity is preserved by the $B_n$ action, so it is enough to construct a vector field of the required form in the neighborhood of each such matrix $S$.

\begin{lemm}\label{68}
Assume that $s_{i,j}=1$ with $i\leq j$, then for each pair 
$(k,\ell)$ in the set
\[
\{(i,j),(i,j+1),\dots , (i,n)\}\cup \{(i+1,j), \dots , 
(j,j),(j,j+1),\dots ,(j,n)\}
\]
there is a vector field $v_{k,\ell}$, vanishing at $S$, 
such that
\begin{enumerate}[(i)]
\item $v_{k,\ell}(\Delta)\in \O\cdot\Delta$, and 
\item the linear part of $v_{k,\ell}$ at $S$ is equal to $(x_{k,\ell}-s_{k,\ell})\frac{\p }{\p {x_{k,\ell}}}$, and in particular is diagonal.	
\end{enumerate}
\end{lemm}

\begin{proof}\
\begin{asparaenum}
\item If $(k,\ell)=(i,\ell)$ with $j< \ell $, then
\[
v_{i,\ell}=x_{i,\ell}\xi_{j,\ell}= x_{i,\ell}\big[ x_{i,j}\p _{i,\ell}
\mod \mm_S\Der \big]=x_{i,\ell}\p _{i,\ell}\mod 
\mm^2_S\Der .
\]

\item If $(k,\ell)=(i,j)$, then if $i<j$
\[
v_{i,j}
=(x_{i,j}-1)\xi_{j,j}=(x_{i,j}-1)\big[ x_{i,j}\p _{i,j}\mod\mm_S\Der \big] 
=(x_{i,j}-1)\p _{i,j}\mod \mm^2_S\Der.
\]
and if $i=j$
\[
v_{i,i}
=\frac{1}{2}(x_{i,i}-1)\xi_{i,i} = (x_{i,i}-1)\p _{i,i}\mod\mm^2_S\Der.
\]  

\item If $(k,\ell)=(k,j)$, with $i< k<j$ then
\[
v_{k,j}
=x_{k,j}\xi_{i,k}= x_{k,j}\big[x_{i,j}\p _{k,j}+\mm_S\Der \big]=x_{k,j}\p _{k,j}\mod\mm^2_S\Der.
\]

\item If $(k,\ell)=(j,j)$ with $i<j$, then
\[
v_{j,j}
=\frac{1}{2}x_{j,j}\xi_{i,j}
=\frac{1}{2}x_{j,j}\big[2 x_{i,j}\p _{j,j}\mod\mm_S\Der \big]
=x_{j,j}\p _{j,j}\mod\mm^2_S\Der.
\]

\item  If $(k,\ell)=(j,\ell)$ with $j< \ell$, then 
\[
v_{j,\ell}=x_{j,\ell}\xi_{i,\ell}
=x_{j,\ell}\big[ x_{i,j}\p _{j,\ell}\mod\mm_S\Der \big]=x_{j,\ell}\p _{j,\ell}\mod\mm^2_S\Der.
\]

\end{asparaenum}
\end{proof}
 
\begin{lemm}\label{69}
For each $i\in\{1,\dots , n\}$ there is a vector field $v_i$ vanishing at 
$S$, such that 
\begin{enumerate}[(i)]
\item $v_i(\Delta)\in\O\cdot\Delta$
\item the linear part of $v_i$ at $S$ is
$
\sum_{k,\ell} \lambda _{k,\ell}(x_{k,\ell}-s_{k,\ell})\frac{\p}{\p {x_{k,\ell}}}
$
where
$\lambda _{k,\ell}=0$ if $k>i$  and 
$\lambda _{i,\ell}> 0$ if $\ell \geq i$; in particular it is
diagonal. 
\end{enumerate}
\end{lemm}

In other words we have a triangular-type system of diagonal linear 
parts with positive terms on the $i$'th row 
and zeros on rows after the $i$'th.

\begin{proof}
If $s_{i,j}=0$ for any $j\geq i$, and $s_{k,i}=0$ for any  
$k\leq i$, we can take $v_i=\xi _{i,i}$.

If $s_{k,i}=1$, with $k\leq i$, then we may apply Lemma~\ref{68} and 
a linear combination of the vector fields $v_{i,i},v_{i,i+1},\dots , 
v_{i,n}$ does the trick.

Finally if $s_{i,j}=1$ for some $j>i$, we observe that $\xi_{i,i}-\xi_{j,j}$, is diagonal and has non zero positive eigenvalues in the positions  
\[
\{(i,i),\dots,(i,j-1)\}\cup\{(i,j+1),\dots(i,n)\}
.\] 
Then we see that the vector field
\[
v_i=v_{i,j}+\xi_{i,i}-\xi_{j,j}+v_{i+1,j}+\cdots + v_{j,j}+v_{j,j+1}+\cdots+v_{j,n}
\]
does the trick since by adding $v_{i,j}$ we complete the row $i$ by a 
positive eigenvalue at $(i,j)$, and we cancel 
with the help of the appropriate $v_{k,\ell}$ all the 
negative eigenvalues with row indices $k>i$.
\end{proof}

\begin{prop}
There is an Euler vector field $v$, $v(\Delta)\in\O\cdot\Delta$ vanishing at $S$, 
with linear part diagonal and having only strictly positive eigenvalues.
\end{prop}

\begin{proof}
We construct $v$, by decreasing induction on $i$, as a linear combination $\alpha _nv_n+ \cdots +\alpha _1v_1$ with positive coefficients, with $\alpha _i>0$ large enough following the choice of $\alpha_n,\dots,\alpha_{i+1}$. 
By construction we have $v(\Delta)=\lambda\Delta$ with $\lambda\in\O$.
\end{proof}

This completes the proof of Proposition~\ref{67}.

\begin{rema}\label{26}
To conclude, we mention a recent theorem of Feh\'er and Patakfalvi. 
In \cite{FP07} they prove that the discriminant $D$ in the representation space of a root of a Dynkin quiver is locally quasihomogeneous. 
Their theorem (\cite[Thm.~5.2]{FP07}) is stated in terms of the Incidence Property that is the subject of their paper, but their proof consists essentially of the contruction of the requisite $\CC^*$-action. 
As a consequence, the LCT holds for these discriminants, by Theorem~\ref{6}.
\end{rema}

\bibliographystyle{cdraifplain}
\bibliography{lfd}

\end{document}